\documentclass[11pt]{amsart}
\usepackage{amsmath}
\usepackage{amssymb}
\usepackage{amsfonts}
\usepackage{amsthm}
\usepackage{enumerate}
\usepackage[mathscr]{eucal}
\usepackage{verbatim}

\begin{document}

\newtheorem{theorem}{Theorem} [section]
\newtheorem{proposition}[theorem]{Proposition}
\newtheorem{conjecture}[theorem]{Conjecture}
\def\theconjecture{\unskip}
\newtheorem{corollary}[theorem]{Corollary}
\newtheorem{lemma}[theorem]{Lemma}
\newtheorem{sublemma}[theorem]{Sublemma}
\newtheorem{observation}[theorem]{Observation}
\theoremstyle{definition}
\newtheorem{definition}{Definition}
\newtheorem{remark}{Remark}
\def\theremark{\unskip}
\newtheorem{question}{Question}
\def\thequestion{\unskip}
\newtheorem{example}{Example}
\def\theexample{\unskip}
\newtheorem{problem}{Problem}

\numberwithin{theorem}{section}
\numberwithin{definition}{section}
\numberwithin{equation}{section}

\theoremstyle{plain}
\newtheorem{thmsub}{Theorem}[subsection]
\newtheorem{lemmasub}[thmsub]{Lemma}
\newtheorem{corollarysub}[thmsub]{Corollary}
\newtheorem{propositionsub}[thmsub]{Proposition}
\newtheorem{defnsub}[thmsub]{Definition}
\numberwithin{equation}{section}

\def\prob{\mu}
\def\fpqs{{F^{pq}_{\sigma}}}
\def\intslash{\rlap{\kern  .32em $\mspace {.5mu}\backslash$ }\int}
\def\qsl{{\rlap{\kern  .32em $\mspace {.5mu}\backslash$ }\int_{Q}}}
\def\Re{\operatorname{Re\,}}
\def\Im{\operatorname{Im\,}}
\def\mx{{\max}}
\def\mn{{\min}}
\def\vth{\vartheta}
\def\tDel{{\widetilde \Delta}}
\def\eann{{{\mathcal E_o}}}
\def\rn{\rr^{n}}
\def\rr{\mathbb R}
\def\R{\mathbb R}
\def\Q{\mathcal Q}
\def\N{\mathbb N}
\def\complex{{\mathbb C}}
\def\norm#1{{ \left|  #1 \right| }}
\def\Norm#1{{ \left\|  #1 \right\| }}
\def\set#1{{ \left\{ #1 \right\} }}
\def\floor#1{{\lfloor #1 \rfloor }}
\def\emph#1{{\it #1 }}
\def\diam{{\text{\rm diam}}}
\def\osc{{\text{\rm osc}}}
\def\ffB{\mathcal B}
\def\itemize#1{\item"{#1}"}
\def\seq{\subseteq}
\def\Id{\text{\sl Id}}

\def\Ga{\Gamma}
\def\ga{\gamma}
\def\Th{\Theta}

\def\prd{{\text{\it prod}}}
\def\parab{{\text{\it parabolic}}}

\def\eg{{\it e.g. }}
\def\cf{{\it cf}}
\def\Rn{{\mathbb R^n}}
\def\Rd{{\mathbb R^d}}
\def\sgn{{\text{\rm sign }}}
\def\rank{{\text{\rm rank }}}
\def\corank{{\text{\rm corank }}}
\def\coker{{\text{\rm coker }}}
\def\loc{{\text{\rm loc}}}
\def\spec{{\text{\rm spec}}}

\def\comp{{\text{\rm comp}}}

\def\Coi{{C^\infty_0}}
\def\dist{{\text{\rm dist}}}
\def\diag{{\text{\rm diag}}}
\def\supp{{\text{\rm supp }}}
\def\rad{{\text{\rm rad}}}
\def\Lip{{\text{\rm Lip}}}
\def\inn#1#2{\langle#1,#2\rangle}
\def\biginn#1#2{\big\langle#1,#2\big\rangle}
\def\rta{\to}
\def\lta{\gets}
\def\noi{\noindent}
\def\lcontr{\rfloor}
\def\lco#1#2{{#1}\lcontr{#2}}
\def\lcoi#1#2{\imath({#1}){#2}}
\def\rco#1#2{{#1}\rcontr{#2}}
\def\bin#1#2{{\pmatrix {#1}\\{#2}\endpmatrix}}
\def\meas{{\text{\rm meas}}}

\def\card{\text{\rm card}}
\def\lc{\lesssim}
\def\gc{\gtrsim}
\def\pv{\text{\rm p.v.}}

\def\alp{\alpha}             \def\Alp{\Alpha}
\def\bet{\beta}
\def\gam{\gamma}             \def\Gam{\Gamma}
\def\del{\delta}             \def\Del{\Delta}
\def\eps{\varepsilon}
\def\ep{\epsilon}
\def\zet{\zeta}
\def\tet{\theta}             \def\Tet{\Theta}
\def\iot{\iota}
\def\kap{\kappa}
\def\ka{\kappa}
\def\lam{\lambda}            \def\Lam{\Lambda}
\def\la{\lambda}             \def\La{\Lambda}
\def\sig{\sigma}             \def\Sig{\Sigma}
\def\si{\sigma}              \def\Si{\Sigma}
\def\vphi{\varphi}
\def\ome{\omega}             \def\Ome{\Omega}
\def\om{\omega}              \def\Om{\Omega}

\def\fA{{\mathfrak {A}}}
\def\fB{{\mathfrak {B}}}
\def\fC{{\mathfrak {C}}}
\def\fD{{\mathfrak {D}}}
\def\fE{{\mathfrak {E}}}
\def\fF{{\mathfrak {F}}}
\def\fG{{\mathfrak {G}}}
\def\fH{{\mathfrak {H}}}
\def\fI{{\mathfrak {I}}}
\def\fJ{{\mathfrak {J}}}
\def\fK{{\mathfrak {K}}}
\def\fL{{\mathfrak {L}}}
\def\fM{{\mathfrak {M}}}
\def\fN{{\mathfrak {N}}}
\def\fO{{\mathfrak {O}}}
\def\fP{{\mathfrak {P}}}
\def\fQ{{\mathfrak {Q}}}
\def\fR{{\mathfrak {R}}}
\def\fS{{\mathfrak {S}}}
\def\fT{{\mathfrak {T}}}
\def\fU{{\mathfrak {U}}}
\def\fV{{\mathfrak {V}}}
\def\fW{{\mathfrak {W}}}
\def\fX{{\mathfrak {X}}}
\def\fY{{\mathfrak {Y}}}
\def\fZ{{\mathfrak {Z}}}

\def\fa{{\mathfrak {a}}}
\def\fb{{\mathfrak {b}}}
\def\fc{{\mathfrak {c}}}
\def\fd{{\mathfrak {d}}}
\def\fe{{\mathfrak {e}}}
\def\ff{{\mathfrak {f}}}
\def\fg{{\mathfrak {g}}}
\def\fh{{\mathfrak {h}}}
\def\fj{{\mathfrak {j}}}
\def\fk{{\mathfrak {k}}}
\def\fl{{\mathfrak {l}}}
\def\fm{{\mathfrak {m}}}
\def\fn{{\mathfrak {n}}}
\def\fo{{\mathfrak {o}}}
\def\fp{{\mathfrak {p}}}
\def\fq{{\mathfrak {q}}}
\def\fr{{\mathfrak {r}}}
\def\fs{{\mathfrak {s}}}
\def\ft{{\mathfrak {t}}}
\def\fu{{\mathfrak {u}}}
\def\fv{{\mathfrak {v}}}
\def\fw{{\mathfrak {w}}}
\def\fx{{\mathfrak {x}}}
\def\fy{{\mathfrak {y}}}
\def\fz{{\mathfrak {z}}}

\def\bbA{{\mathbb {A}}}
\def\bbB{{\mathbb {B}}}
\def\bbC{{\mathbb {C}}}
\def\bbD{{\mathbb {D}}}
\def\bbE{{\mathbb {E}}}
\def\bbF{{\mathbb {F}}}
\def\bbG{{\mathbb {G}}}
\def\bbH{{\mathbb {H}}}
\def\bbI{{\mathbb {I}}}
\def\bbJ{{\mathbb {J}}}
\def\bbK{{\mathbb {K}}}
\def\bbL{{\mathbb {L}}}
\def\bbM{{\mathbb {M}}}
\def\bbN{{\mathbb {N}}}
\def\bbO{{\mathbb {O}}}
\def\bbP{{\mathbb {P}}}
\def\bbQ{{\mathbb {Q}}}
\def\bbR{{\mathbb {R}}}
\def\bbS{{\mathbb {S}}}
\def\bbT{{\mathbb {T}}}
\def\bbU{{\mathbb {U}}}
\def\bbV{{\mathbb {V}}}
\def\bbW{{\mathbb {W}}}
\def\bbX{{\mathbb {X}}}
\def\bbY{{\mathbb {Y}}}
\def\bbZ{{\mathbb {Z}}}

\def\cA{{\mathcal {A}}}
\def\cB{{\mathcal {B}}}
\def\cC{{\mathcal {C}}}
\def\cD{{\mathcal {D}}}
\def\cE{{\mathcal {E}}}
\def\cF{{\mathcal {F}}}
\def\cG{{\mathcal {G}}}
\def\cH{{\mathcal {H}}}
\def\cI{{\mathcal {I}}}
\def\cJ{{\mathcal {J}}}
\def\cK{{\mathcal {K}}}
\def\cL{{\mathcal {L}}}
\def\cM{{\mathcal {M}}}
\def\cN{{\mathcal {N}}}
\def\cO{{\mathcal {O}}}
\def\cP{{\mathcal {P}}}
\def\cQ{{\mathcal {Q}}}
\def\cR{{\mathcal {R}}}
\def\cS{{\mathcal {S}}}
\def\cT{{\mathcal {T}}}
\def\cU{{\mathcal {U}}}
\def\cV{{\mathcal {V}}}
\def\cW{{\mathcal {W}}}
\def\cX{{\mathcal {X}}}
\def\cY{{\mathcal {Y}}}
\def\cZ{{\mathcal {Z}}}

\def\tA{{\widetilde{A}}}
\def\tB{{\widetilde{B}}}
\def\tC{{\widetilde{C}}}
\def\tD{{\widetilde{D}}}
\def\tE{{\widetilde{E}}}
\def\tF{{\widetilde{F}}}
\def\tG{{\widetilde{G}}}
\def\tH{{\widetilde{H}}}
\def\tI{{\widetilde{I}}}
\def\tJ{{\widetilde{J}}}
\def\tK{{\widetilde{K}}}
\def\tL{{\widetilde{L}}}
\def\tM{{\widetilde{M}}}
\def\tN{{\widetilde{N}}}
\def\tO{{\widetilde{O}}}
\def\tP{{\widetilde{P}}}
\def\tQ{{\widetilde{Q}}}
\def\tR{{\widetilde{R}}}
\def\tS{{\widetilde{S}}}
\def\tT{{\widetilde{T}}}
\def\tU{{\widetilde{U}}}
\def\tV{{\widetilde{V}}}
\def\tW{{\widetilde{W}}}
\def\tX{{\widetilde{X}}}
\def\tY{{\widetilde{Y}}}
\def\tZ{{\widetilde{Z}}}

\def\tcA{{\widetilde{\mathcal {A}}}}
\def\tcB{{\widetilde{\mathcal {B}}}}
\def\tcC{{\widetilde{\mathcal {C}}}}
\def\tcD{{\widetilde{\mathcal {D}}}}
\def\tcE{{\widetilde{\mathcal {E}}}}
\def\tcF{{\widetilde{\mathcal {F}}}}
\def\tcG{{\widetilde{\mathcal {G}}}}
\def\tcH{{\widetilde{\mathcal {H}}}}
\def\tcI{{\widetilde{\mathcal {I}}}}
\def\tcJ{{\widetilde{\mathcal {J}}}}
\def\tcK{{\widetilde{\mathcal {K}}}}
\def\tcL{{\widetilde{\mathcal {L}}}}
\def\tcM{{\widetilde{\mathcal {M}}}}
\def\tcN{{\widetilde{\mathcal {N}}}}
\def\tcO{{\widetilde{\mathcal {O}}}}
\def\tcP{{\widetilde{\mathcal {P}}}}
\def\tcQ{{\widetilde{\mathcal {Q}}}}
\def\tcR{{\widetilde{\mathcal {R}}}}
\def\tcS{{\widetilde{\mathcal {S}}}}
\def\tcT{{\widetilde{\mathcal {T}}}}
\def\tcU{{\widetilde{\mathcal {U}}}}
\def\tcV{{\widetilde{\mathcal {V}}}}
\def\tcW{{\widetilde{\mathcal {W}}}}
\def\tcX{{\widetilde{\mathcal {X}}}}
\def\tcY{{\widetilde{\mathcal {Y}}}}
\def\tcZ{{\widetilde{\mathcal {Z}}}}

\def\tfA{{\widetilde{\mathfrak {A}}}}
\def\tfB{{\widetilde{\mathfrak {B}}}}
\def\tfC{{\widetilde{\mathfrak {C}}}}
\def\tfD{{\widetilde{\mathfrak {D}}}}
\def\tfE{{\widetilde{\mathfrak {E}}}}
\def\tfF{{\widetilde{\mathfrak {F}}}}
\def\tfG{{\widetilde{\mathfrak {G}}}}
\def\tfH{{\widetilde{\mathfrak {H}}}}
\def\tfI{{\widetilde{\mathfrak {I}}}}
\def\tfJ{{\widetilde{\mathfrak {J}}}}
\def\tfK{{\widetilde{\mathfrak {K}}}}
\def\tfL{{\widetilde{\mathfrak {L}}}}
\def\tfM{{\widetilde{\mathfrak {M}}}}
\def\tfN{{\widetilde{\mathfrak {N}}}}
\def\tfO{{\widetilde{\mathfrak {O}}}}
\def\tfP{{\widetilde{\mathfrak {P}}}}
\def\tfQ{{\widetilde{\mathfrak {Q}}}}
\def\tfR{{\widetilde{\mathfrak {R}}}}
\def\tfS{{\widetilde{\mathfrak {S}}}}
\def\tfT{{\widetilde{\mathfrak {T}}}}
\def\tfU{{\widetilde{\mathfrak {U}}}}
\def\tfV{{\widetilde{\mathfrak {V}}}}
\def\tfW{{\widetilde{\mathfrak {W}}}}
\def\tfX{{\widetilde{\mathfrak {X}}}}
\def\tfY{{\widetilde{\mathfrak {Y}}}}
\def\tfZ{{\widetilde{\mathfrak {Z}}}}


\def\ta{{\tilde a}}
\def\tb{{\tilde b}}
\def\tc{{\tilde c}}
\def\tx{{\tilde x}}
\def\ty{{\tilde y}}
\def\tz{{\tilde z}}
\def\tnu{{\tilde \nu}}
\def\tn{{\tilde n}}
\def\tf{{\tilde f}}
\def\ts{{\tilde s}}
\def\tt{{\tilde t}}
\def\tsi{{\tilde \sigma}}
\def\tsigma{{\tilde \sigma}}

\def\Atil{{\widetilde A}}          \def\atil{{\tilde a}}
\def\Btil{{\widetilde B}}          \def\btil{{\tilde b}}
\def\Ctil{{\widetilde C}}          \def\ctil{{\tilde c}}
\def\Dtil{{\widetilde D}}          \def\dtil{{\tilde d}}
\def\Etil{{\widetilde E}}          \def\etil{{\tilde e}}
\def\Ftil{{\widetilde F}}          \def\ftil{{\tilde f}}
\def\Gtil{{\widetilde G}}          \def\gtil{{\tilde g}}
\def\Htil{{\widetilde H}}          \def\htil{{\tilde h}}
\def\Itil{{\widetilde I}}          \def\itil{{\tilde i}}
\def\Jtil{{\widetilde J}}          \def\jtil{{\tilde j}}
\def\Ktil{{\widetilde K}}          \def\ktil{{\tilde k}}
\def\Ltil{{\widetilde L}}          \def\ltil{{\tilde l}}
\def\Mtil{{\widetilde M}}          \def\mtil{{\tilde m}}
\def\Ntil{{\widetilde N}}          \def\ntil{{\tilde n}}
\def\Otil{{\widetilde O}}          \def\otil{{\tilde o}}
\def\Ptil{{\widetilde P}}          \def\ptil{{\tilde p}}
\def\Qtil{{\widetilde Q}}          \def\qtil{{\tilde q}}
\def\Rtil{{\widetilde R}}          \def\rtil{{\tilde r}}
\def\Stil{{\widetilde S}}          \def\stil{{\tilde s}}
\def\Ttil{{\widetilde T}}          \def\ttil{{\tilde t}}
\def\Util{{\widetilde U}}          \def\util{{\tilde u}}
\def\Vtil{{\widetilde V}}          \def\vtil{{\tilde v}}
\def\Wtil{{\widetilde W}}          \def\wtil{{\tilde w}}
\def\Xtil{{\widetilde X}}          \def\xtil{{\tilde x}}
\def\Ytil{{\widetilde Y}}          \def\ytil{{\tilde y}}
\def\Ztil{{\widetilde Z}}          \def\ztil{{\tilde z}}


\def\ahat{{\hat a}}          \def\Ahat{{\widehat A}}
\def\bhat{{\hat b}}          \def\Bhat{{\widehat B}}
\def\chat{{\hat c}}          \def\Chat{{\widehat C}}
\def\dhat{{\hat d}}          \def\Dhat{{\widehat D}}
\def\ehat{{\hat e}}          \def\Ehat{{\widehat E}}
\def\fhat{{\hat f}}          \def\Fhat{{\widehat F}}

\def\ghat{{\hat g}}          \def\Ghat{{\widehat G}}
\def\hhat{{\hat h}}          \def\Hhat{{\widehat H}}
\def\ihat{{\hat i}}          \def\Ihat{{\widehat I}}
\def\jhat{{\hat j}}          \def\Jhat{{\widehat J}}
\def\khat{{\hat k}}          \def\Khat{{\widehat K}}
\def\lhat{{\hat l}}          \def\Lhat{{\widehat L}}
\def\mhat{{\hat m}}          \def\Mhat{{\widehat M}}
\def\nhat{{\hat n}}          \def\Nhat{{\widehat N}}
\def\ohat{{\hat o}}          \def\Ohat{{\widehat O}}
\def\phat{{\hat p}}          \def\Phat{{\widehat P}}
\def\qhat{{\hat q}}          \def\Qhat{{\widehat Q}}
\def\rhat{{\hat r}}          \def\Rhat{{\widehat R}}
\def\shat{{\hat s}}          \def\Shat{{\widehat S}}
\def\that{{\hat t}}          \def\That{{\widehat T}}
\def\uhat{{\hat u}}          \def\Uhat{{\widehat U}}
\def\vhat{{\hat v}}          \def\Vhat{{\widehat V}}
\def\what{{\hat w}}          \def\What{{\widehat W}}
\def\xhat{{\hat x}}          \def\Xhat{{\widehat X}}
\def\yhat{{\hat y}}          \def\Yhat{{\widehat Y}}
\def\zhat{{\hat z}}          \def\Zhat{{\widehat Z}}


\def\intprod{\mathbin{\lr54}}
\def\reals{{\mathbb R}}
\def\integers{{\mathbb Z}}
\def\complex{{\mathbb C}\/}
\def\naturals{{\mathbb N}\/}
\def\distance{\operatorname{distance}\,}
\def\degree{\operatorname{degree}\,}
\def\dim{\operatorname{dimension}\,}
\def\Span{\operatorname{span}\,}
\def\ZZ{ {\mathbb Z} }
\def\e{\varepsilon}
\def\p{\partial}
\def\rp{{ ^{-1} }}
\def\Re{\operatorname{Re\,} }
\def\Im{\operatorname{Im\,} }
\def\ov{\overline}
\def\bx{{\bf{x}}}
\def\lt{L^2}
\def\Farrow{F} 
\def\Phiarrow{\Phi} 
\def\expect{\mathbb E}

\def\scriptx{{\mathcal X}}
\def\scriptb{{\mathcal B}}
\def\scripta{{\mathcal A}}
\def\scriptk{{\mathcal K}}
\def\scriptd{{\mathcal D}}
\def\scriptp{{\mathcal P}}
\def\scriptl{{\mathcal L}}
\def\scriptv{{\mathcal V}}
\def\scripti{{\mathcal I}}
\def\scripth{{\mathcal H}}
\def\scriptm{{\mathcal M}}
\def\scripte{{\mathcal E}}
\def\scriptt{{\mathcal T}}
\def\scriptb{{\mathcal B}}
\def\frakg{{\mathfrak g}}
\def\frakG{{\mathfrak G}}

\author{Jonathan Bennett}
\address{Jonathan Bennett\\
School of Mathematics\\University of Birmingham\\
The Watson Building,
Edgbaston\\
Birmingham, B15 2TT, United Kingdom} \email{J.Bennett@bham.ac.uk}

\author{Andreas Seeger}
\address{
        Andreas Seeger\\
        Department of Mathematics\\
        University of Wisconsin\\
    Madison, Wisconsin 53706-1388, USA}
\email{seeger@math.wisc.edu}

\thanks{Research partially supported by EPSRC Postdoctoral Fellowship
GR/S27009/02   (J.B.) and
by  the National Science Foundation (A.S.)}


\title[Fourier extension on spheres  and  oscillatory integrals]
{The Fourier extension operator   on large spheres \\  and related oscillatory
integrals}


\begin{abstract}
We obtain new estimates for a class of oscillatory integral operators
with folding canonical relations satisfying a curvature condition. The main
lower
bounds showing sharpness are proved using Kakeya set constructions. As a
special case of the upper bounds we deduce optimal
$L^p(\bbS^2)\to L^q(R\bbS^2)$ estimates
for the Fourier extension operator
on large spheres in $\mathbb{R}^3$, which are uniform in the radius $R$.
Two appendices are included, one  concerning an application
to Lorentz space bounds for
averaging operators along curves in $\bbR^3$,
and one on  bilinear estimates.
\end{abstract}

\maketitle

\section{Introduction} \label{intro}
For  functions $g\in L^1(\mathbb{S}^{d})$ on the $d$-dimensional
unit sphere we define the Fourier extension operator to be the mapping
$\cE: g\mapsto \widehat {gd\sigma} $ where
$$
 \widehat {gd\sigma}
(\xi)=\int_{\mathbb{S}^d}e^{- i\inn{x}{\xi}}g(x)d\sigma(x),
$$
$d\sigma$ denotes the rotation invariant measure
on $\mathbb{S}^d$
induced by  Lebesgue measure in $\bbR^{d+1}$,
 and $\xi\in\mathbb{R}^{d+1}$. We note that
the adjoint of this operator is the Fourier
restriction operator
$f\mapsto\widehat{f}\bigl|_{\mathbb{S}^d}$, where $\;\widehat{\;}\;$ denotes
the
Euclidean Fourier transform in  $d+1$ dimensions.
A substantial amount of recent work is concerned with weighted inequalities
of the general
form
\begin{equation}\label{genweightineq}
\Big(\int |\widehat {gd\sigma}|^q d\mu\Big)^{1/q}
\lc \|g\|_{L^p(\mathbb{S}^d)}
\end{equation}
for certain measures $\mu$ on
$\bbR^{d+1}$. \footnote{Throughout this paper we will use the notation
$X\lesssim Y$ ($X\gtrsim Y$) if for non-negative quantities $X$ and $Y$
there exists a constant $C>0$ such that $X\leq CY$ ($X\geq CY$). The dependence
of the implicit constant $C$ on various parameters present will be clarified
by the context.}
Perhaps the most notable instance of this is the case of Lebesgue measure,
which corresponds to the
classical Fourier restriction problem; see for example
\cite{Feff}, \cite{stein},  \cite{Tomas}, \cite{Bo-Bes} and
\cite{Tao}. In addition to this, the inequalities
\eqref{genweightineq} for
certain broader classes of measures $\mu$ are known to have
applications to a variety of
well-known and largely unsolved
problems
in partial differential equations,
harmonic analysis
and geometric measure theory;
see
\cite{BRV}, \cite{RV}, \cite{CS1}, \cite{CS2},
\cite{Wo}, \cite{carse}, \cite{SjSo}, \cite{Ma}, \cite{ErRev}, \cite{Er},
and many
further references contained in those papers.
The content of the current paper is partially motivated by the particular situation where
the measures $\mu$ are supported on
large spheres in $\mathbb{R}^{d+1}$; this has been
studied recently in \cite{BBC1},
\cite{BBC2} and \cite{BCSV}. We take $\mu$ to be the rotation invariant measure
on $\la S^{d}$ induced by Lebesgue measure in $\bbR^{d+1}$.
In particular  the case for  circles in the plane is well understood;
namely the
$L^p(\bbS^1)\to L^q (\la \mathbb{S}^1)$ operator norm of $\cE$ is
uniformly  bounded in $\la$, if,  and only if,  $q\ge 3$ and
$p\ge q/(q-2)$. This follows from a result on more general
oscillatory integral operators  in \cite{GS3};  for further discussion and an
alternative proof of the $L^3$ bound see \cite{BCSV}.
Here we prove for spheres in $\bbR^3$:

\begin{theorem}\label{2Dextthm}
The inequality
\begin{equation}\label{2Dext}
\big\|\widehat{g d\sigma}\big\|_{L^q(\la \bbS^2)} \le C \|g\|_{L^p(\bbS^2)}
\end{equation} holds for all $\la$,  all $g\in L^p(\bbS^2)$ and some $C$, if
and only if $q>5/2$ and $p\ge 2q/(2q-3)$.
\end{theorem}


After rescaling   one sees that
uniform $L^p(\bbS^d)\to L^q(\la \bbS^d)$ bounds for $\cE$
are equivalent with  the $O(\la^{-d/q})$ bound for the
$L^p(\bbS^d)\to L^q(\bbS^d)$ operator norm of $\cE_\la$, given by
\begin{equation}\label{ERdef}\cE_\la g(\xi)=\widehat {gd\sigma}(\la\xi).
\end{equation}

%
%
The operators  $\cE$ and $\cE_\la$ are closely related to a
 Radon transform arising in scattering theory, considered by
 Melrose and Taylor \cite{MT}.
After appropriately parametrizing $\mathbb{S}^d$
the operator $\cE_\la$
may be seen as a special case of a much
more general class of oscillatory  operators
acting on functions defined on $\bbR^d$, given by
%
%
\begin{equation}\label{generaloscop}
T_\la f(x)=\int e^{i\la\phi(x,y)} \chi(x,y) f(y) dy.
\end{equation}
Here $\phi$ is a smooth real-valued phase function on
$\Omega_L\times\Omega_R$ where $\Omega_L$ and $\Omega_R$ are open
subsets of $\mathbb{R}^d$ and $\chi$ is smooth with
$\supp\chi\subset\Omega_L\times\Omega_R$. We shall now discuss the assumptions
on  the phase  which are appropriate for the study of $\cE_\la$.

The $L^2$ mapping properties of $T_\la$ are governed by geometrical properties
of the  canonical relation associated to the phase
$\phi$; it  is defined to be the
(twisted) graph of the gradient map,
$$\cC_\phi=\{(x,\nabla_x\phi,y,-\nabla_y\phi):(x,y)\in\supp\chi\}\subset T^*\Omega_L\times T^*\Omega_R.$$
Here we assume that the
 projections $\pi_L$ and $\pi_R$ mapping $\cC_\phi$
to $ T^*\Omega_L$ and $T^*\Om_R$,  respectively,
\begin{equation}
\begin{aligned}
&\pi_L: (x,y)\mapsto (x,\phi_x(x,y))
\\
&\pi_R: (x,y)\mapsto (x,\phi_y(x,y))
\end{aligned}
\end{equation}
are Whitney folds.
Analytically the fold condition on $\pi_L$ can be expressed by
requiring that $\corank d\pi_{L}\le 1$  and when $\dim \ker d\pi_L=1$
then the Hessian considered as a map from $\ker d\pi_L$ to $\coker d\pi_L$ is nonzero; {\it i.e.}
\begin{equation} \label{piLfoldhess}
\ 0\neq b\in \ker \phi_{xy},
\ 0\neq a\in \coker \phi_{xy},
  \quad \implies\quad \inn{b}{\nabla_y}^2 \inn{a}{\phi_x}\neq 0 .
\end{equation}
An equivalent condition is
\begin{equation} \label{piLfold}
\det \phi_{xy} (x,y)=0, \ 0\neq b\in \ker \phi_{xy},
  \quad \implies\quad \inn{b}{\nabla_y} (\det\phi_{xy})\neq 0.
\end{equation}

Similarly the corresponding condition on $\pi_R$ being a Whitney fold is
\begin{equation} \label{piRfold}
\det \phi_{yx} (x,y)=0, \ 0\neq a\in \ker \phi_{yx},
  \quad \implies\quad \inn{a}{\nabla_x} (\det\phi_{yx})\neq 0.
\end{equation}
Using the terminology in \cite{MT}
we say that $\cC_\phi$
 is a
folding canonical relation if \eqref{piLfold} and
\eqref{piRfold} are satisfied.
The $L^2$ operator norm of $T_\la$ is $O(\la^{-d/2+1/6})$ by the work of
Melrose-Taylor \cite{MT} and Pan-Sogge \cite{PanSo}.


Condition \eqref{piLfold} makes
\begin{equation} \label{defofL}\cL=\{(x,y): \det \phi_{xy}=0\}
\end{equation}
a smooth hypersurface in
$\bbR^d\times \bbR^{d}$; moreover for fixed $x$
\begin{equation}\label{Lfreezingx}
\{y: (x,y)\in \cL\}\end{equation}
is a smooth hypersurface in $\bbR^d$, and thus the varieties
\begin{equation}
\cL_x:= \{\xi \in T^*_x \Omega_L: \xi= \phi_x(x, y), \ (x,y)\in \cL\}
\end{equation}
are smooth hypersurfaces in the fibers.
Following \cite{GS1} we assume the following condition
(which is based on the Carleson-Sj\"olin hypothesis,
\cf. \cite{H}, \cite{MSS}):

\medskip

\noi{\bf Curvature condition:}
\begin{equation}
\label {curvcond}
\text{\emph{ For every  $x\in \Omega_L$, the
hypersurface $ \cL_x$
is  convex and has nonvanishing curvature.}}
\end{equation}
The convexity  and nonvanishing curvature hypotheses
 mean that the second fundamental form is either positive definite
or negative definite everywhere on $\cL_x$.

Condition \eqref{curvcond}  is not relevant for $L^2\to L^2$
bounds, however it is crucial for $L^p\to L^q$ bounds in higher dimensions.
In one dimension there is no curvature condition
and the best possible results are known, namely
\begin{equation}\label{oneD}
\big\|T_\la\|_{L^p(\bbR)\to L^q(\bbR)} \lc\la^{-1/q},
\quad
q\geq 2p', \, q\geq 3,
\end{equation}
holds under the assumptions \eqref{piLfold}, \eqref{piRfold}. This was proved in \cite{GS3}.
Examples
(see \S\ref{lowerbds}) show that the sharp bound
\begin{equation}\label{dbound}
\big\|T_\la\|_{L^p(\bbR^d)\to L^q(\bbR^d)} \lc\la^{-d/q},
\end{equation}
can only hold for $q\geq (d+1)p'/d$   and $q\geq(2d+1)/d$
(here $p'=p/(p-1)$). In two and higher dimensions
Kakeya type examples  exclude the case $q=(2d+1)/d$.
Under assumption \eqref{curvcond} inequality \eqref{dbound}  has been
established  by Greenleaf and one of the  authors \cite{GS1} in the
range $q\geq (2d+2)/d$;
actually in \cite{GS1} the assumption of a
folding canonical relation has been replaced with a weaker
 one-sided assumption involving only the projection $\pi_L$.
Moreover,  in the range
$q\geq (2d+2)/d$ the definiteness assumption on the second
fundamental form can be replaced by merely the nondegeneracy assumption
(of course this makes no difference when $d=2$).

Under the folding relation and curvature  assumptions
we improve the known range $q\ge 3$ of inequality  \eqref{dbound} in two dimensions,
and get a best possible result.

\begin{theorem} \label{main2d}
Suppose that $d=2$ that $\cC_\phi$ is a folding canonical relation
 and  that the curvature condition \eqref{curvcond}  is satisfied.
Then for $\la\geq 2$
\begin{equation}\label{qgr52}
\big\|T_\la\|_{L^p(\bbR^2)\to L^q(\bbR^2)} \lc\la^{-2/q}, \quad
q\ge  \frac{3p'}{2}, \quad q>\frac 52.
\end{equation}

Moreover,
\begin{equation}\label{ple52}
\big\|T_\la\|_{L^{q}(\bbR^2)\to L^{q}(\bbR^2)}\lc
 \la^{-\frac 23-\frac{1}{3q}} (\log\la )^{\frac 12-\frac 1q}, \quad 2\leq q<5/2,
\end{equation}
and
\begin{equation}\label{rwt}
\big\|T_\la\|_{L^{5/2,1}(\bbR^2)\to L^{5/2,\infty}(\bbR^2)}\lc
 \la^{-4/5} (\log\la )^{1/10}.
\end{equation}

The estimates are sharp in the following sense: If  there is a
point $P\in \cL$ so that
$\chi(P)\neq 0$ then there is a positive constant $c>0$ depending on $\chi$  and $\la_0>1$
so that for all  $\la\ge \la_0$
\begin{equation}\label{sharpness}
\big\|T_\la\|_{L^{p,1}(\bbR^2)\to L^{q,\infty}(\bbR^2)}\ge c
\max\{\la^{-\frac 2q},
 \la^{-\frac{2}{3p'}-\frac 1q},\,
\la^{-\frac23-\frac{1}{3q}}  (\log  \la)^{\frac 12-\frac 1q} \,\}.
\end{equation}
\end{theorem}

It would
 be interesting to know whether the restricted weak type estimate \eqref{rwt}
could  be replaced by an
 $L^{5/2}\to L^{5/2}$ estimate with the same bounds; this
remains open.

The  assumptions of Theorem
\ref{main2d} are satisfied for the operator $\cE_\la$ in the Fourier extension problem on spheres, so that Theorem
\ref{2Dextthm} is a direct consequence of Theorem \ref{main2d}
(see \S\ref{preparations}).
Indeed the spheres on both sides of inequality
\eqref{2Dext} may be replaced by compact pieces of two surfaces  in $\bbR^3$
with nonvanishing Gaussian curvature.


\medskip

{\it Structure of the paper.} In \S\ref{preparations} we discuss
some preparatory   changes of variables which are useful in the
proof of both the necessary and sufficient conditions,
 and briefly
discuss the validity of our assumptions for the phases in the
Fourier extension problem. In
\S\ref{lowerbds} we prove the sharpness of Theorem \ref{main2d};
the main part of this section is concerned with a Kakeya type
example. In \S\ref{basicdecomposition} we give the basic
decompositions of the operator in terms of the size of
$\det\phi_{xy}''$ and state the main estimates for these pieces.
In \S\ref{models} we discuss easy proofs of the required bounds in
certain model cases and raise some open questions. The more
technical proof of the main estimates
 in the general case is given in \S\ref{proofpartone}, \S\ref{proofoffirstprop}
and \S\ref{B0bound}.
The paper has two appendices. In the first one, \S\ref{appI}, we consider the
convolution with measures on some curves in $\bbR^3$;
we use a variant of our estimates to give a Lorentz-space improvement of
 Oberlin's endpoint estimates \cite{Ob}. In the second appendix, \S\ref{appII},
 we revisit the bilinear
estimates from \cite{BBC1} and give a straightforward proof based on the
geometric properties of the canonical relation.

\medskip

\noi{\bf Acknowledgements:}

J.B.  would like to thank Juan Antonio Barcel\'o, Tony Carbery,
Fernando Soria and Ana Vargas for their on-going collaborative
work on the subject of general weighted $L^2$ norm inequalities
for the Fourier extension operator. The $L^p(\mathbb{S}^d)\to
L^q(R\mathbb{S}^d)$ Fourier extension problems addressed in this
paper arose naturally in this work. Thanks are also due to Tony
Carbery for his involvement in the early stages of this project.

A.S. would like to thank Allan Greenleaf for numerous  conversations on
oscillatory integral operators, many of them related to this project,  in
the course of their long term  collaboration.

\

\section{Preparation of the phase function}\label{preparations}

It is advantageous to suitably
prepare the phase function by possibly changing variables
in $x$ and in $y$.
 These changes of variables affect  the estimates
only by constants. We have to observe that our hypotheses
are invariant  under these changes  of variables.
This is standard for the conditions \eqref{piLfold} and
\eqref{piRfold}.
Concerning the curvature condition
a change of  variables in $x$  induces a linear change in the fiber
($\xi$-)~variables and thus leaves the curvature condition invariant.
We now examine  the independence of parametrization and
 invariance under change of the $y$-variable, of the curvature condition.
We shall consider the situation in $d$ dimensions.

If $x$ is fixed and $z\mapsto G(x,z)$ is a regular parametrization
of $\{y:(x,y)\in \cL\}$  (with parameter $z\in \bbR^{d-1}$) then vectors in
$\coker \phi_{xy}$ are normal  to the hypersurface $\cL_x$ in the fiber
above $x$  and
 the curvature condition
is just saying that for $v\in\coker \phi_{xy}$
the Hessian of the
map
$$z\mapsto \inn{v}{\nabla_x\phi(x, G(x,z))}$$
is either positive definite or negative definite;
 this Hessian equals
$$\tfrac{\partial G}{\partial z}^T
\inn{v}{\nabla_x\phi}''_{yy}
\tfrac{\partial G}{\partial z}+ O(\inn{v}{\phi_{xy}})
$$
at $y=G(x,z)$ and the last term drops out since
$v\in\coker\phi_{xy} (x,G(x,z))$. From this the invariance easily follows.



We now prepare our phase function to have an approximate
normal form at a point $P=(x^o,y^o)$, and we may assume that
$P\in\cL$. (\textit{i.e. } to have certain  derivatives vanish at $P$). Let
us assume that the phase function  $$(x,y)\mapsto \psi(x,y)$$ has
a canonical relation  $\cC_\psi$ satisfying \eqref{piLfold},
\eqref{piRfold}, and \eqref{curvcond}. We shall find
diffeomorphisms $G_L$ and $G_R$, mapping neighborhoods  of the
origins of $\bbR^d_L$, $\bbR^d_R$, to neighborhoods of $x^o$,
$y^o$ respectively, so that at the origin  $O=(O_L,O_R)$ the phase
\begin{equation}\label{phifrompsi}
\phi(x,y)=\psi(G_L(x),G_R(y)) \end{equation}
satisfies
the conditions
\begin{equation}
\det \phi_{x'y'}(O)\neq 0
\label{detprime}
\end{equation}
and
\begin{align}
&\phi_{xy_d}(O)= 0
\label{phixyd}
\\
&\phi_{x_dy}(O)= 0;
\label{phixdy}
\end{align}
moreover
\begin{align}
&\phi_{x_dy_dy_d}(O)\neq  0,
\label{phixdydyd}
\\
&\phi_{x_dx_dy_d}(O)\neq  0,
\label{phixdxdyd}
\\
&\phi_{x_d y_d y'}(O)=  0,
\label{mixedydterm}
\\
&\phi_{x_dy_d x'}(O)=  0,
\label{mixedxdterm}
\end{align}
and also
\begin{equation}\label{straightnormal}
\phi_{x'y'x_d}(O)=0.
\end{equation}

To accomplish this, let $a$ and $b$ be unit vectors in $\bbR^d_L$ and $\bbR^d_R$ respectively, so that at $P$ we have
$\psi_{xy}b=0$, $a^T\psi_{xy}=0$
(recall that we assume that at $P$ the kernel and cokernel of
$\psi_{xy}$ are one dimensional). Now choose rotations
$\rho_L$ of $\bbR^d_L$ and
$\rho_R$ of $\bbR^d_R$  so that $\rho_L^{-1}b=e_d$,
$\rho_R^{-1}a=e_d$. Then
$$\phi^{[1]}(x,y)= \psi(x^o+ \rho_L(x),y^o+\rho_R(y))$$ satisfies
$\det \phi^{[1]}_{x'y'}(O)\neq 0$,
$\phi^{[1]}_{x'y_d}(O)= 0$ and $\phi^{[1]}_{x_dy'}(O)= 0$.
By the formula
\begin{equation}\label{splitdet}
\det \phi_{xy}= \det(\phi_{x'y'}) (\phi_{x_dy_d}- \phi_{x_dy'}\phi_{x'y'}^{-1}\phi_{x'y_d})
\end{equation}
(applied to $\phi^{[1]}$)
we also have
$\phi^{[1]}_{x_dy_d}(O)= 0$
and see that $\phi^{[1]} $ satisfies
\eqref{detprime}, \eqref{phixyd} and \eqref{phixdy}.
Notice that from the fold assumptions \eqref{piLfold} and
\eqref{piRfold} we also have
$\phi^{[1]}_{x_dy_dy_d}(O)\neq  0
$, and $\phi^{[1]}_{x_dx_dy_d}(O)\neq  0$.

We now consider the phase-function
$$\phi^{[2]}(x,y)=\phi^{[1]}(\sigma_L (x), \sigma_R(y))$$
for suitable shears in
$\bbR^d_L$ and $\bbR^d_R$, of the form
\begin{align*}
\sigma_L(x)=(x', x_d-\sum_{i=1}^{d-1} \alpha_i x_i), \notag \qquad
\sigma_R(y)= (y', y_d-\sum_{j=1}^{d-1} \beta_j y_j). \notag
\end{align*}
Note that $\phi^{[2]}$ still satisfies
\eqref{detprime}, \eqref{phixyd} and \eqref{phixdy}, and also
\eqref{phixdydyd} and \eqref{phixdxdyd}, independently
of the choice of $\alpha$ and $\beta$.
Now if  we choose
\begin{equation*}
\alpha_i=-\frac{\phi^{[1]}_{x_dy_d x_i}(O)}{\phi^{[1]}_{x_dy_dx_d}(O)},
\quad\beta_j=-\frac{\phi^{[1]}_{x_dy_d y_j}(O)}
{\phi^{[1]}_{x_dy_dy_d}(O)},
\end{equation*}
then conditions
\eqref{mixedydterm}, \eqref{mixedxdterm} are satisfied for $\phi^{[2]}$
 as well.

Now set
\begin{equation*}
\phi(x,y)=\phi^{[2]}((x'+x_d Bx', x_d),y)
\end{equation*}
where $B=-[\phi^{[2]}_{x'y'}(O)]^{-1}\phi^{[2]}_{y'x'x_d}(O)$.
Then
\eqref{phifrompsi} holds with
$G_R(y) =y^0+ \rho_R(\sigma_R(y))$ and
$G_L(x) =x^0+ \rho_L(\sigma_L(\nu(x)))$,  where
$\nu(x)=
(x'+x_d Bx', x_d)$. The phase
$\phi$ satisfies  \eqref{straightnormal} and
conditions \eqref{detprime} --
\eqref{mixedxdterm} continue to hold.

Finally, by replacing the phase
 $\phi(x,y)$ with
$\phi(x,y)-\phi(x^o,y)$
 we may assume that
\begin{equation}
\label{yderiv}
\partial^\alpha_y \phi (O)=0,
\end{equation}
for all multiindices $\alpha$.

We now examine the curvature condition \eqref{curvcond} at $O$.
By  condition \eqref{phixdydyd}  we can
solve near $O$
\begin{equation}
\label{implicit}
\det \phi_{xy}=0 \quad \iff  \quad y_d=g(x,y')
\end{equation}
with $g(O_L,O_R')=0$.
Implicit differentiation and condition
\eqref{mixedydterm} implies that
\begin{equation}
\label{gradg}
\nabla_{y'} g(O_L, O_R')=0.
\end{equation}
Thus our curvature condition at $O$ reads
\begin{equation}
\label{curvatP}
\nabla_{y'y'}^2 \big(\phi_{x_d}(O_L, y', g(O_L,y')\big)\Big|_{y'=O_R'}
\text{ is positive or negative   definite,}
\end{equation}
which by \eqref{mixedydterm} and \eqref{gradg} reduces to the
definiteness
assumption on the Hessian of $\phi_{x_d}$, namely,
\begin{equation}
\label{curvatPreduced}
\nabla^2_{y'y'}\phi_{x_d}(P)
\text{ is positive or negative   definite.}
\end{equation}

\subsubsection*{On the phase functions in the  Fourier extension problem}
We briefly discuss here how the extension operator $\mathcal{E}_{\lambda}$
in \eqref{ERdef} of the introduction belongs to our general family of
oscillatory integral operators $T_{\lambda}$ satisfying \eqref{piLfold},
\eqref{piRfold} and \eqref{curvcond}.

Let $S$ be a patch of a smooth
convex hypersurface of
$\bbR^{d+1}$, with nonvanishing Gaussian curvature
(in particular $S$ may be part of $\bbS^d$ as in \eqref{ERdef}).
Let $y\mapsto \Gamma(y)$ be a
parametrization of $S$ (where the parameter $y$ is chosen from an
open subset of $\bbR^d$). Let $\Sigma$ be a smooth hypersurface of
$\Bbb R^{d+1}$, parametrized by $x\mapsto \Xi(x)$, where $x$ belongs to
an open set of $\bbR^{d}$.
Then the operator $\cE_\la$ in \eqref{ERdef} may now be written as an
oscillatory integral operator with phase function
\begin{equation} \label{extphase}
\phi(x,y)= \inn{\Xi(x)}{ \Gamma(y)}.
\end{equation}
Clearly $\phi_{xy}=\Xi'(x)^T \Gamma'(y)$ is of rank $\ge d-1$ and
$\cL$ consists of those
$(x,y)$ for which
the normal line
for $S$
at $\Gamma(y)$ is parallel to the tangent space for $\Sigma$ at $\Xi(x)$
(or, equivalently,
the normal line
for $\Sigma$
at $\Xi(x)$ is parallel to the tangent space for $S$ at $\Gamma(y)$).

The assumption that the second fundamental form of $S$ is definite
implies that the fold condition for $\pi_L$, \eqref{piLfold},
is satisfied.
Indeed if
$(x,y)\in \cL$ and
if $a$, $b$ are
nonzero vectors in $\bbR^d$ so that
$a^T\Xi'(x)^T
\Gamma'(y)=0$ and
$\Xi'(x)^T \Gamma'(y)b=0$ then the fold condition in the form
 \eqref{piLfoldhess} is
saying that
$$\inn{b}{\nabla_y} \,a^T\Xi'(x)^T\Gamma'(y)b\neq 0$$
and this is implied by the definiteness of the fundamental form of $S$
 since
$\Xi'(x) a$ is a nonzero vector
 perpendicular to the tangent space of $S$ at $\Gamma(y)$.

For the curvature condition \eqref{curvcond} we fix $x$ and
solve
$\det \phi_{xy}=0$  by $y=G(x,z)$
so that
$\cL_x$ is parametrized by $z\mapsto \Xi'(x)^T\Gamma(G(x,z))$.
We need to verify that  the second fundamental form of
$\cL_x$ is definite, {\it i.e.} that
$$\nabla_{zz}^2 \inn{a}{\Xi'(x)^T\Gamma(G(x,z))}$$
is definite if $\inn{a}{\Xi'(x)^T\Gamma'(G(x,z))}=0$.
However under this last  condition the second fundamental form becomes
$$\tfrac{\partial G}{\partial z}^T
  \nabla_{yy}^2 \inn{\Xi'(x)a}{\Gamma(y)}
\tfrac{\partial G}{\partial z} \quad\text{at } y=G(x,z).$$
Again as $\Xi'(x) a$ is normal to $S$ at $\Gamma(y)$
we see by the definiteness
assumption on the second fundamental form and
by $\rank\frac{\partial G}{\partial z}=d-1$,
that the last displayed formula gives a
 definite  $(d-1)\times(d-1)$ matrix. Thus the curvature condition is verified.

Finally, if in addition
 we also assume that $\Sigma$ is a convex hypersurface
with   nonvanishing curvature then we see by symmetry that
the fold condition for $\pi_R$, \eqref{piRfold}, is satisfied as well
(see
 also  \cite{com-rad} for a discussion of the  structure of
$\pi_R$ in the  more general
situation where $\Sigma$ is convex and of finite line type).

\section{Lower bounds} \label{lowerbds}
We now establish  lower bounds for  the operator norms of $T_\la$
showing in particular the sharpness of Theorem \ref{main2d}. We
work in $d$ dimensions and assume that the fold and curvature
conditions \eqref{piLfold}, \eqref{piRfold} and \eqref{curvcond}
hold, and in addition we make the  (necessary)  assumption  that
there is a point
 $(x^o,y^o)\in \cL$
for which
\begin{equation}\label{chinonzero}\chi(x^o,y^o)\neq 0.
\end{equation}

By the reductions described
in \S \ref{preparations} we may assume that $(x^o,y^o)=O$, that
(\ref{detprime}-\ref{straightnormal}) hold,
and, in dimension $d\ge 2$, that
\eqref{curvatPreduced} holds.

%

We are interested in the range of exponents $(p,q)$ for which
\begin{equation} \label{ladq}
\|T_\la \|_{L^{p,1}
\to L^{q,\infty}}
\lc \la^{-d/q}
\end{equation} holds.
It is easy to see that the decay rate in  \eqref{ladq} is sharp
(for any $C^1$  phase function). Since the operator is local we
have
\begin{equation}\label{locality}
\|T_\la \|_{L^{\infty}
\to L^{q,\infty}}\lc
\|T_\la \|_{L^{p,1}
\to L^{q,\infty}}
\end{equation}
and therefore it suffices to prove lower
bounds for the weak type $(\infty,q)$ operator norm. Without loss
of generality
 $\Re(\chi(x,y))>c>0$ for $|x|\le \eps$, $|y|\le \eps$. Let $\la\gg \eps^{-1}$, and define $f(y)=e^{-i\la\phi(0,y)}$ for $|y|\le \eps$ and $f(y)=0$ elsewhere.
Then $|T_\la f(x)|\ge c>0$ for $|x|\le c_0\eps\la^{-1}$ and thus
\begin{equation}
\label{dqlwbd}
\|T_\la\|_{L^\infty\to L^{q,\infty}}\ge c'\la^{-d/q}.
\end{equation}

 The following simple
lemma shows that the condition $q\le (d+1)p'/d$ is necessary   for
\eqref{ladq} to hold. Note that \eqref{dqlwbd} and
\eqref{simplestlow} yield the first two lower bounds stated in
\eqref{sharpness}.

\begin{lemma}\label{onetube}
There is $c>0$ so that
\begin{equation}
\label{simplestlow}
\|T_\la\|_{L^{p,1}\to L^{q,\infty}}\ge c
\la^{d/(3p)-d/3-(2d-1)/(3q)}.
\end{equation}
\end{lemma}

\begin{proof}
 Let $f_0$ be
the characteristic function of the ball $\{y:|y|\leq \eps
\la^{-1/3}\}$, and define
$$f(y)= f_0(y)\exp(-i\la(\inn{y}{\phi_y(O)}+ \tfrac 12\inn{y}{\phi_{yy}(O)y}))$$
so that $\|f\|_{L^{p,1}} \approx \la^{-d/(3p)}$. By considering
the Taylor expansion of $\phi(x,y)-\phi(x,0)$ we observe that
$$|\phi(x,y)-\phi(x,0)-\langle y,\phi_y(O)\rangle-\tfrac{1}{2}\langle y,\phi_{yy}(O)y\rangle|\le C\eps\la^{-1}$$
whenever $|y|\le \eps \la^{-1/3}$, $|x|\le \eps\la^{-1/3} $ and
$|\inn {x}{\phi_{xy}(O)y}|\le \la^{-1}$. On multiplying
$T_{\lambda}f(x)$ by the unimodular factor
$e^{-i\lambda\phi(x,0)}$, we find that if $x$ is such that these
conditions hold uniformly in $|y|\leq\eps\la^{-1/3}$, then
$$|T_\la f(x)|\geq c\lambda^{-d/3}$$
if $\eps$
 is sufficiently small.
By the assumptions \eqref{phixyd} and \eqref{phixdy} we see that
$|\inn {x}{\phi_{xy}(O)y}|\le \la^{-1}$ holds for all
$|y|\leq\eps\la^{-1/3}$ whenever $|x'|\leq \eps' \la^{-2/3}$ and
$|x_d|\leq \eps' \lambda^{-1/3}$. Thus $\|T_\la
f\|_{L^{q,\infty}}\gc \la^{-d/3-(2d-1)/(3q)}$, and the assertion
follows.
\end{proof}

We shall now show by a  randomization argument that for $d\geq 2$
the inequality \eqref{ladq} can only hold for $q> (2d+1)/d$, and also
establish the sharpness of \eqref{ple52}, \eqref{rwt}. The approach is
inspired
by the result of Beckner, Carbery, Semmes and Soria  \cite{BCSS}
on the failure of
restricted weak type endpoint bounds for the classical Fourier extension
operator
({\it cf.}  also Tao's generalization  \cite{Tao-CS}
to oscillatory integral operators).
We
use a rescaled version of  the Kakeya construction in Keich
\cite{Ke}. Let $\delta\ll 1$, $\delta<\alpha< 1/10$ and suppose
that   for every $n'\in \bbZ^{d-1}$ with $|n'\delta|\le \alpha$ we
are given a $r\delta\times\dots\times r\delta\times r$ rectangle
$P_n$ passing through the
 hyperplane
$x_d=0$  so that the long edges are parallel to $(n'\delta,1)$.
Then there are vectors $v_n\in \bbR^{d-1}\times\{0\}$, $|v_n|\le
|\alpha|$,
 so that the union of translated rectangles  $v_n+P_n$
satisfies
\begin{equation}\label{compression}
\big|\bigcup_n  v_n+P_n\big|\le C \big(\log (\alpha/\delta)\big)^{-1}
\sum_{n}|P_n|.
\end{equation}
We shall apply this fact after possible changes of variables,
 with $r=\delta=\la^{-1/3}$, $\alpha=\epsilon  \lambda^{-1/6}$,
and large $\la$, then $\log (\alpha/\delta)\approx \log \la$.

\begin{proposition} \label{kh}
Suppose  $d\ge 2$ and $q>2$, then  there is
$c>0$ and $\la_0>0$ so that  for all $\lambda>\lambda_0$
\begin{equation}
\label{kkbd}
\|T_\la\|_{L^{\infty}\to L^{q,\infty}}\ge c
\la^{-d/3-(d-1)/(3q)}  (\log  \la)^{1/2-1/q}.
\end{equation}
\end{proposition}

\begin{proof}
Assume without loss of generality $\Re (\chi(x,y))\ge 1$ whenever
$|x_i|\le c_0$, $|y_j|\le c_0$ for some constant $c_0>0$. Assume
$\la\gg c_0^{-1}$. We let  $\fQ$ be the family of all cubes of
$\mathbb{R}^d$ of sidelength $\la^{-1/6}$, of the form
$\prod_{i=1}^d[n_i\la^{-1/6}, (n_i+1)\la^{-1/6})$ where
$n=(n_1,\hdots,n_d)\in\mathbb{R}^d$ and $|n_i|\la^{-1/6} \le
c_0/2$. For $Q\in \fQ$ let $x_Q$ be the center of $Q$. Let
$y_Q=(x_Q', g(x_Q, x_Q'))$ where $g$ is given by \eqref{implicit},
and let $B(Q)$ be the ball of radius $\eps_1 \la^{-1/6}$ centered
at $y_Q$. In view of \eqref{phixdydyd} and \eqref{phixdxdyd} we
have $g_{x_d}\neq 0$ near the origin and by choosing $c_0$
sufficiently small we may assume that $y\mapsto (y', g(y,y'))$ is
a diffeomorphism near the origin. Consequently, if $\eps_1$ is
sufficiently small, the balls $\{B(Q): Q\in \fQ\}$ form a disjoint
family.

On each cube $Q$, and each ball $B(Q)$
we shall now change variables as in \S\ref{preparations}.
Namely, for each $Q\in\fQ$ there is  a diffeomorphism $\fv_Q$
mapping a neighborhood
$\cU_{L,Q}$
of the origin $O_L$ to an open set  $\cV_Q$ containing $Q$
and a diffeomorphism
$\fw_Q$
mapping a neighborhood $\cU_{R,Q}$ of the origin $O_R$ to a neighborhood
 $\cW_Q$   of $y_Q$ containing $B(Q)$,
so that
the phase function
$$
\psi^Q(x,y)=\phi(\fv_Q(x),\fw_Q(y))
$$
satisfies  conditions (\ref{detprime}-\ref{straightnormal}), and
also \eqref{curvatPreduced} holds for $\psi^Q$. The  bounds for
the derivatives of $\psi^Q$ are uniform in $Q$, as are the
implicit lower bounds in \eqref{detprime}, \eqref{phixdydyd},
\eqref{phixdxdyd}, \eqref{curvatPreduced}  for those functions. We
can find a  positive $\eps_2 \ll \eps_1$ so that for every $Q$
the sets $\cU_{L,Q}$
and $\cU_{R,Q}$ contain the cubes of sidelength $\eps_2$  centered
at the origins $O_L$ and $O_R$, respectively.
Moreover there is a positive $\eps_3
\le \eps_2$, so that if $Q_o$ denotes the cube of sidelength
$\eps_3\la^{-1/6}$ centered at $O_L$
 then
$\fv_Q(Q_o)\subset Q$, for every $Q\in \fQ$.
We let
$\cZ=\{n\in \bbZ^d: |n|\le 10^{-1}\eps_3 \la^{1/6}\}$.

We decompose  this cube $Q_o$  into plates at height $\la^{-1/3}
n_d$, with $|n_d|\le \eps_3 10^{-1} \la^{1/6}$. Let $\Pi_{n'}$ be the
orthogonal projection to the hyperplane orthogonal to $(n'
\la^{-1/3},1)$. We now apply the above mentioned construction by
Keich (with angular parameter $\alpha\le\eps_4 \la^{-1/6}$, {\it
cf.} \eqref{compression}). Then for each $n_d$ we find a family of
$\lambda^{-2/3}\times\cdots\times\la^{-2/3}\times\la^{-1/3}$
rectangles $\widetilde R_{n}=\widetilde{R}_{n',n_d}$
so that $\widetilde R_n$
contains the set
$$
R_n=\{x:|x_d-n_d\la^{-1/3}|\le \la^{-1/3}\eps_4,\,
|\Pi_{n'}
(x-a(n))|\le \eps_4 \la^{-2/3}\}
$$
where
\begin{equation}\label{defofan}
a(n)=(a'(n),a_d(n)) \in \Bbb R^{d-1}\times \{\la^{-1/3}n_d\},
\text{  with } |a'(n)|\le \eps_3\la^{-1/6},
\end{equation}
 and, for the measure of
$$E(n_d)=
\bigcup_{n'} R_{n', n_d},$$
there is the Besicovich type estimate
\begin{equation}\label{measkakeya}
|E(n_d)|\le
C \frac{\la^{(d-1)/6}\la^{-(2d-1)/3}}{\log \la},
\end{equation}
uniformly in $n_d$. Observe that for $n\in\mathbb{Z}^d$ the
rectangle $\widetilde{R}_{n}$ lies in the plate at height
$n'\lambda^{-1/3}$, contains the point $a(n)\in\mathbb{R}^d$ and
has long sides in the direction $(n'\lambda^{-1/3},1)$.




\begin{sublemma}\label{Tlafnlem}
If $\eps\ll \eps_4$ is sufficiently small then there is $c(\eps)
>0$ so that
 the following holds for  $\la\ge \eps^{-1}$.
For each $Q\in \fQ$ there is a disjoint family of balls $B_{n,Q}$,
$n\in \cZ$, each of radius $\eps \la^{-1/3}$ and contained in
$\fw_Q^{-1}(B(Q))$, and for each $(Q,n)\in \fQ\times \cZ$ there is
a smooth function $H_{n,Q}$ defined on $B_{n,Q}$ so that with
$$f_{n,Q}(y)=\chi_{B_{n,Q}}(\fw_Q^{-1}(y)) e^{-i\la H_{n,Q}(\fw_Q^{-1}(y))}$$
we have
\begin{equation}\label{Tlafn}
|T_{\la} f_{n,Q}(x)|\ge c(\eps) \la^{-d/3}, \quad\text{ if }
x\in \cR_{n,Q}: = \fv_Q( R_n).
\end{equation}
\end{sublemma}

We postpone the proof of the sublemma and continue with the proof of
the proposition.
We show
\begin{equation}\label{Khconsequence}
\Bigl\|\sum_{n,Q}\chi_{\cR_{n,Q}}
\Bigr\|_{L^{q/2,\infty}(\mathbb{R}^d)}
\lesssim\lambda^{2d/3-d/(3p)} \|T_\la\|_{L^\infty\to L^{q,\infty}}^2.
\end{equation}
To see \eqref{Khconsequence}
we follow the argument in \cite{BCSS}.
We denote by $\{r_k\}$ the system of Rademacher functions.
Choose an injective function
$(n,Q) \mapsto k(n,Q)$ with values in the positive integers.
By Khinchine's inequality
$$
\Big(\sum_{n,Q}|T_\la f_{n,Q}(x)|^2\Big)^{1/2}\lc
\int_0^1\Big|\sum_{n,Q} r_{k(n,Q)}(t) T_\la f_{n,Q}(x) \Big| dt
$$ uniformly in $x$, $\la$.
Now by the sublemma
$\chi_{\cR_{n,Q}} \lc \la^{2d/3}|T_\la f_{n,Q}(x)|^2$,
and hence
\begin{align*}
&\Bigl\|\sum_{n,Q}\chi_{\cR_{n,Q}}
\Bigr\|_{L^{q/2,\infty}(\mathbb{R}^d)}\lc
\la^{2d/3}\Big\|\sum_{n,Q}|T_\la f_{n,Q}|^2\Big\|_{L^{q/2,\infty}}
\\&=
\la^{2d/3}\Big\|\Big(\sum_{n,Q}|T_\la f_{n,Q}|^2\Big)^{1/2}\Big\|_{L^{q,\infty}}^2
   \lc
\la^{2d/3}\Big\|\int_0^1\Big|
\sum_{n,Q} r_{k(n,Q)}(t)T_\la f_{n,Q}\Big| dt\Big\|_{L^{q,\infty}}^2
\end{align*}
and the square root of the right hand side is further estimated by
a constant times
\begin{align*}
\la^{d/3}\Big\|\int_0^1\Big|T_\la&\big[
\sum_n r_{k(n,Q)}(t) f_{n,Q}\big]\Big| dt\Big\|_{L^{q,\infty}}
\le \la^{d/3}
\int_0^1\Big\|
T_\la[\sum_{n,Q} r_{k(n,Q)}(t) f_{n,Q}]\Big\|_{L^{q,\infty}} dt
\\
&\le \la^{d/3}\|T_\la\|_{L^\infty\to L^{q,\infty}}
\int_0^1\Big\|
\sum_{n,Q} r_{k(n,Q)}(t) f_{n,Q}\Big\|_{L^\infty} dt
\\&\lc \la^{d/3} \|T_\la\|_{L^\infty\to L^{q,\infty}}
\Big\|
\sum_{n,Q} |f_{n,Q}|\Big\|_{L^{\infty}}
\\&\lc\la^{d/3}
\|T_\la\|_{L^{\infty}\to L^{q,\infty}}.
\end{align*}
For the last inequality we have used the disjointness of the supports of
$f_{n,Q}$ which follows from the disjointness of the balls $B(Q)$, and for each fixed $Q$ from the disjointness of the $B_{n,Q}$, $n\in \cZ$.

Next observe  that $|\cR_{n,Q}|\approx |R_n|\approx \la^{-2(d-1)/3-1/3}$
and that $\card(\cZ)\approx \la^{d/6}$, $\card(\fQ)\approx \la^{d/6}$ and
therefore
$$ \la^{-(d-1)/3}
\approx
\Bigl\|\sum_{Q\in \fQ}\sum_{n\in \cZ} \chi_{\cR_{n,Q}}\Bigr\|_1.$$
Hence,
by  the duality of $L^{(q/2)',1}$ and $L^{q/2,\infty}$
and \eqref{measkakeya},
\begin{eqnarray*}
\begin{aligned}
\la^{-(d-1)/3}
&\lesssim
\Big|\bigcup_{Q\in\fQ}\bigcup_{(n',n_d)\in\cZ} \fv_Q(R_n)\Big|^{1-2/q}
\Bigl\|\sum_{Q\in \fQ}\sum_{n\in \cZ} \chi_{\cR_{n,Q}}\Bigr\|_{L^{q/2,\infty}}
\\
&\lc\Big[ \sum_{Q\in \fQ}\sum_{|n_d|\lc \la^{1/6}}\big|E(n_d)|\Big]^{1-2/q}
\la^{2d/3}\|T_\la\|_{L^\infty\to L^{q,\infty}}
\\
&\lc \Big(\la^{(d+1)/6} \frac{\la^{(d-1)/6}\la^{-(2d-1)/3}}{\log \la}
\Big)^{1-2/q}\la^{2d/3}\|T_\la\|_{L^\infty\to L^{q,\infty}}^2,
\end{aligned}
\end{eqnarray*}
which implies
$$\|T_\la\|_{L^{\infty}\to L^{q,\infty}}\ge c \la^{-d/3-(d-1)/(3q)}
(\log \la)^{1/2-1/q}
$$
and thus the assertion.
\end{proof}

\begin{proof}[Proof of Sublemma \ref{Tlafnlem}]
We fix $Q$; our estimates will be uniform in $Q$ and  we will
generally suppress indices indicating the dependence of the terms
on $Q$. For $f$ defined near $O_R$ (in particular   in $B(Q)$) and
for  $x\in Q$ we set $$\cT_\la f(x)= T_\la[f\circ
\fw_Q^{-1}](\fv_Q(x)).$$ Then
$$
\cT_\la f(x)=\int e^{i\la \psi(x,y)}\chi_1(x,y) f(y) dy,$$
where
$\chi_1(x,y)=\chi(\fv_Q(x), \fw_Q(y))|\det \fw_Q'(y)| $,  and the phase
$$\psi(x,y)\equiv \psi^Q(x,y)=\phi(\fv_Q(x), \fw_Q(y))$$
satisfies conditions (\ref{detprime}-\ref{straightnormal}).
We also note that $\det \psi(x,y)=0$ when $y_d=\fg(x,y')$
and $\fg$ satisfies  \eqref{gradg} and \eqref{curvatP}
(with $\phi$ replaced by $\psi$).

We shall now identify balls $B_n$ so that for suitable $f_n$
supported on $B_n$ the function  $\cT_\la f_n$ is bounded below by
$c\la^{-d/3}$ on $R_n$. To achieve this we argue very much as in
the proof of Lemma \ref{onetube} and analyze the Taylor expansion
of the phase function $\psi(x,y)-\psi(x,b)$ about $x=a$, for
suitable $a$, $b$. Let
\begin{equation}
\label{defofH}
H(a,b,y)=\psi_y(a,b) (y-b) +\frac 12 (y-b)^t\psi_{yy}(a,b)(y-b).
\end{equation}
Then
\begin{equation}\label{taylorpsiexp}
\psi(x,y)-\psi(x,b)=
H(a,b,y)+ (x-a)^t\psi_{xy}(a,b)(y-b) +
O(|y-b||x-a|^2+|y-b|^3))
\end{equation}
and we further split with $\psi^{y'x'}:=\psi_{x'y'}^{-1}$
\begin{equation}
\label{taylorexpansion}
\begin{aligned}
(x-a)^t\psi_{xy}(y-b) =&
\big(
(x'-a')^t+(x_d-a_d) \psi_{x_dy'}\psi^{y'x'}  \big)
\big(
\psi_{x'y'}(y'-b')+\psi_{x'y_d}(y_d-b_d)
\big)
\\
&+ (x_d-a_d) \big(\psi_{x_dy_d}- \psi_{x_dy'}\psi^{y'x'}\psi_{x'y_d}
\big)(y_d-b_d)
\end{aligned}
\end{equation}
where the derivatives of $\psi$ are evaluated at $(a,b)$.
Note that by \eqref{splitdet} the second term drops out if
$b_d=\fg(a,b')$.

To define $B_n$ and $f_n$ we first consider for fixed $a_d$ the map
$\sigma(\cdot, a_d)$ defined in a neighborhood of the origin
$O_R'$ of $\bbR^{d-1}$ by
$$y'\mapsto \sigma(y',a_d):=
-\psi_{x_d y'}\psi^{y'x'}\big|_{(x',x_d,y',y_d)=
(O_L', a_d, y', \fg(O_L',a_d,y'))}.
$$
Then $\sigma(O_R',0)=O_R'$.
By the curvature condition \eqref{curvatP}
and \eqref{phixdy},
the map $\sigma(\cdot,a_d)$  is a
diffeomorphism on a neighborhood of $O_R'$,
if $a_d$ is small;
the  bounds are  uniform for  $a_d$ in an open interval
containing $0$.
We may assume that the neighborhood of $O_R'$ and its image contain the ball of
radius $\eps_3$ centered at $O_R'$, whenever $|a_d|\le \eps_3$.
Let $b'(n)$ be defined by
$$\sigma(b'(n),\la^{-1/3}n_d)= \la^{-1/3} n'$$
and we assume that $|n|\lc \eps\la^{1/6}$.
Let $$b(n)= (b'(n), b_d(n)):= (b'(n), g(O_L', n_d\la^{-1/3},b'(n))
$$  and let $B_n$ be the ball of radius $\eps\la^{-1/3} $
centered at $b(n)$.
Define
$$f_n(y)=\chi_{B_n}(y) e^{-i\la H(a(n), b(n),y)}$$
with $H$ as in \eqref{defofH}.
It will be crucial to note that
$|\det\psi_{xy}(x,y)|\lc \la^{-1/3}$ when $y\in B_n$, $x\in R_n$
(see \eqref{smalldet} below).

It now suffices to show
\begin{equation}\label{cTlafn}
|\cT_{\la} f_n(x)|\ge c \la^{-d/3}, \quad\text{ if } x\in R_n.
\end{equation}
with  the positive constant $c$ independent of $\la$, $Q$ and $n$.
To see \eqref{cTlafn} note that
$$e^{-i\la \phi(x,b(n))}T_\la f_n(x) =\int e^{i\la \Psi_n(x,y)}
\chi_1(x,y) \chi_{B_n}(y) dy$$
where, by \eqref {taylorpsiexp},
$
\Psi_n(x,y)=(x-a)^t\psi_{xy}(a,b)(y-b)+
O(\eps \la^{-1})$
evaluated at $(a,b)=(a(n),b(n))$, and the error bounds hold if
$|y-a(n)|\le \eps \la^{-1/3}$, and $ |x-a(n)|\le
\la^{-1/3}$.
Thus estimate
\eqref{cTlafn} follows if we verify that
\begin{equation}\label{quadrform}
\big|(x-a(n))^t\psi_{xy}(a(n),b(n))(y-b(n))\big|\le
C\eps \la^{-1}, \text{ if } x\in R_n.
\end{equation}
Since the vector
$$(\la^{-1/3} n',1)= -\psi_{x_dy'}\psi^{y'x'}(0',\la^{-1/3}n_d, b'(n),
 \fg(0,\la^{-1/3}n_d, b'(n)))
$$ is  in the kernel of the orthogonal projection $\Pi_{n'}$
we have for $x\in R_n$
$$|x'-a'(n)+(x_d-a_d(n))
\psi_{x_dy'}\psi^{y'x'}(0',\la^{-1/3}n_d, b'(n),
\fg(0,\la^{-1/3}n_d, b'(n)))|
\le C\eps\la^{-2/3}.
$$
Notice that by the crucial properties \eqref{straightnormal} and
\eqref{phixdy}, \eqref{mixedydterm}  the terms
$\psi_{x_dy'x'}$,
$\psi_{x_d y'}$ and $\psi_{x_dy'y_d}$ are all $O(\eps\la^{-1/6})$ in $Q_o$.
Thus
$$
|\psi_{x_dy'}\psi^{y'x'} (a', a_d, b', \fg(a,b')) -
\psi_{x_dy'}\psi^{y'x'} (0', a_d, b',  \fg(0,a_d,b'))| \le
C\eps\la^{-1/3}.
$$
Consequently
\begin{equation}
\label{quadrform1}
\big|\big((x'-a')^t+(x_d-a_d) \psi_{x_dy'}\psi^{y'x'}(a,b)  \big)
\big(
\psi_{x'y'}(y'-b')+\psi_{x'y_d}(a,b)(y_d-b_d)
\big)\big|\le C\eps\la^{-1}
\end{equation}
if $a=a(n)$, $b=b(n)$, the derivatives are evaluated at $(a(n), b(n))$ and
$x\in R_n$ and $y\in B_n$.
Moreover for these choices of $a,b,x,y$
\begin{equation}
\label{quadrform2}
\big|(x_d-a_d) \big(\psi_{x_dy_d}- \psi_{x_dy'}\psi^{y'x'}\psi_{x'y_d}
\big)(y_d-b_d)\big| \le C\eps\la^{-1}.
\end{equation}
To see this we use that
$b_d(n)-\fg(a(n), b'(n))= \inn{\fg_{x'}(O_L',a_d,b'(n))}{a'(n)} +
O(\la^{-1/3})$ and
since by implicit differentiation using \eqref{mixedxdterm} we have
$\fg_{x'}=O(\la^{-1/6})$ we see that in fact
$$b_d(n)-\fg(a(n), b'(n))=O(\la^{-1/3}).$$
Hence
\begin{equation}\label{smalldet}\psi_{x_dy_d}- \psi_{x_dy'}\psi^{y'x'}\psi_{x'y_d}
=O(\la^{-1/3})\end{equation}
and thus \eqref{quadrform2} follows.
By \eqref{quadrform1} and \eqref{quadrform2} we get \eqref{quadrform} and
this finishes the verification of
\eqref{cTlafn}.
\end{proof}

\subsubsection*{Remark}
The reader familiar with the  wave packet analysis in the
context of the classical restriction problem for the Fourier
transform (see for example \cite{Tao}) may find it enlightening to
construct Kakeya set examples of this type for the particular
operator $g\mapsto\widehat{gd\sigma}\bigl|_{R\mathbb{S}^d}$
discussed in the introduction.
The key point here is that if $B$ is an
$R^{-1/3}$-cap centered at a point $x_0$ on the equator of $\mathbb{S}^d$,
and if $\nu\in\mathbb{S}^{d}$ lies within
a distance of $O(R^{-1/3})$ of the north pole, then the
function $g(x)=\chi_B(x)e^{iR\nu\cdot x}$ is such that
$|\widehat{gd\sigma}\bigl|_{R\mathbb{S}^d}(\xi)|\gtrsim
R^{-d/3}\chi_{T_B}(\xi)$;
here $T_B$ is an ``eccentric cap'' (or ``stretched cap'')
on $R\mathbb{S}^d$ of dimensions $O(R^{1/3})\times\cdots\times
O(R^{1/3})\times O(R^{2/3})$, centered at $R\nu$
and with long edges in the
direction $x_0$.
In order to exploit this we
let $\{\nu_m\}_{1\leq m\lesssim R^{1/3}}$
be a sequence of equally spaced points on the curve
$\{u=(u_1,\hdots,u_{d+1})\in\mathbb{S}^{d}:
u_1=\cdots=u_{d-1}=0\},$ and let $\mathbb{S}_{\nu_m}^{d-1}=
\{\omega\in\mathbb{S}^d:\omega\cdot\nu_m=0\}$. We now choose a collection
of disjoint $R^{-1/3}$-caps $\{B_{m,n}\}_{1\leq m\lesssim R^{1/3},
1\leq n\lesssim R^{(d-1)/3}}$ on $\mathbb{S}^d$ such that for each
$m$ and $n$ the center of $B_{m,n}$ (which we will call $x_{m,n}$) lies on
the great sphere $\mathbb{S}_{\nu_m}^{d-1}$. Now, for each $m$ and $n$ let
$T_{m,n}$ denote an eccentric cap on $R\mathbb{S}^d$
of dimensions $O(R^{1/3})\times\cdots\times O(R^{1/3})\times O(R^{2/3})$,
with long sides pointing in the direction $x_{m,n}$ and centered at a
point $R\nu_{m,n}\in R\mathbb{S}^d$ with $|\nu_{m,n}-\nu_m|\lesssim
R^{-1/3}$. Now if $g_{m,n}(x)=e^{iR\nu_{m,n}\cdot x}
\chi_{B_{m,n}}(x)$ then
$|\widehat{g_{m,n}d\sigma}\bigl|_{R\mathbb{S}^d}(\xi)|\gtrsim
R^{-d/3}\chi_{T_{m,n}}(\xi)$ uniformly in $m$ and $n$. Choosing the
caps $B_{m,n}$ and frequencies $\nu_{m,n}$ appropriately,
taking $g$ to be a random
combination of the form $\sum\pm g_{m,n}$ and invoking
appropriate Besicovitch type estimates now leads to the required
necessary condition $q>(2d+1)/d$. Here of course the $O(R^{1/3})$ scaled
Kakeya sets that feature are subsets of $R\mathbb{S}^d$ rather
than $\mathbb{R}^d$. Notice also that an analogue of the additional
decomposition at scale $O(\lambda^{-1/6})$, required in the treatment
of the general operators $T_{\lambda}$, is not necessary here.
\medskip

\section{Basic decompositions
}\label{basicdecomposition}

It is standard to decompose the operator $T_\la$ in terms of the
size of $\det \phi_{xy}$.
By a Taylor expansion (using \eqref{implicit},
\eqref{phixdydyd})
we observe that  on the (small) support  of our cutoff function
\begin{equation}
y_d-g(x,y')= \cC(x,y)\det\phi_{xy}
\end{equation}
with $\cC(x,y)\neq 0$
so that the decomposition in terms of
$\det\phi_{xy}$ can be realized by decomposing in terms of the size of
$y_d-g(x,y')$.
Thus we split
$T_\la= \sum_{2^l< \la ^{1/3}}
T_{\la,l} + \widetilde T_\la$
where
\begin{equation} \label{Tlaldef}T_{\la,l} f(x)= \int e^{i\la \phi(x,y)} \chi(x,y)
\chi_1(2^l(y_d-g(x,y'))) f(y) dy
\end{equation}
where $\chi_1$ is supported in $(2/3, 3/2)\cup((-3/2,-2/3)$ and
$\widetilde T_\la$
is defined similarly with a cutoff
$\chi_0(2^l(y_d-g(x,y'))$
localizing to the region $|y_d-g(x,y')|\lc \la^{-1/3}$.
Then $T_{\la,l}$
and $\widetilde T_\la$ cover the situations where
$|\det\phi_{xy}|\approx 2^{-l}$,  and
$|\det\phi_{xy}|\lc \la^{-1/3}$, respectively.

By standard $L^2$ theory  \cite{Cu}, \cite{GS3} (see also \cite{PS} for earlier results in special cases)  we have
\begin{align}
\|T_{\la,l}\|_{L^2\to L^2} &\lc 2^{l/2} \la ^{-d/2}, \quad 2^l<\la^{1/3},
\label{L2Tlestimate}
\\
\|\widetilde T_\la\|_{L^2\to L^2} &\lc  \la ^{-(d-1)/2-1/3}.
\label{L2Tlestimatelim}
\end{align}

Our main estimates in two dimensions are
\begin{align}
\|T_{\la,l}\|_{L^4(\bbR^2)\to L^4(\bbR^2)} &\lc 2^{-3l/4} \la ^{-1/2}
 (\log \la)^{1/4},
\quad 2^l<\la^{1/3},
\label{Tlestimate}
\\
\|\widetilde T_\la\|_{L^4(\bbR^2)\to L^4(\bbR^2)} &\lc  \la ^{-3/4}
(\log\la)^{1/4},
\label{Tlestimatelim}
\end{align}
and for $2\le p<4$, $ q=3p'$,
\begin{align}
\|T_{\la,l}\|_{L^p(\bbR^2)\to L^q(\bbR^2)} &\lc 2^{-l/p'} \la ^{-2/q},
\quad 2^l<\la^{1/3},
\label{Tlpqestimate}
\\
\|\widetilde T_\la\|_{L^p(\bbR^2)\to L^q(\bbR^2)} &\lc  \la ^{-1/(3p')-2/q}.
\label{Tlpqestimatelim}
\end{align}

Notice that \eqref{Tlestimatelim},
\eqref{Tlpqestimatelim}
are  limiting cases of
\eqref{Tlestimate} and \eqref{Tlpqestimate}.
We shall prove only
\eqref{Tlestimate} and \eqref{Tlpqestimate}
and the proofs of \eqref{Tlestimatelim}
and \eqref{Tlpqestimatelim} are analogous.
Indeed for the proofs of \eqref{Tlestimatelim}, \eqref{Tlpqestimatelim}
the localization to
the region where $|y_2-g(x,y_1)|\approx 2^{-l}$ can be replaced by
the localization to
the region where $|y_2-g(x,y_1)|\lc 2^{-l}$.

By interpolation it follows from
\eqref{L2Tlestimatelim} (with $d=2$), and
\eqref{Tlestimatelim}
that

\begin{align}
\label{tildeT52estimate}
&\|\widetilde T_\la\|_{L^{5/2}(\bbR^2)\to L^{5/2}(\bbR^2)} \lc  \la ^{-4/5}
(\log \la)^{1/10},
\\
\label{tildeTpqestimate}
&\|\widetilde T_\la\|_{L^{p}(\bbR^2)\to L^{q}(\bbR^2)} \lc  \la ^{-2/q} ,
\quad q=\frac{3p'}2,\quad  q>\frac 52.
\end{align}
Moreover the restricted weak type estimates
\begin{align}
\Big\|\sum_{2^l\leq \la^{1/3}}
T_{\la,l}\Big\|_{L^{5/2,1}(\bbR^2)\to L^{5/2,\infty}(\bbR^2)} \lc  \la ^{-4/5}
(\log \la)^{1/10},
\label{rwtone}
\\
\Big\|\sum_{2^l\leq \la^{1/3}}
T_{\la,l}\Big\|_{L^{p,1}(\bbR^2)\to L^{q,\infty}(\bbR^2)} \lc  \la ^{-2/q} ,\quad q=\frac{3p'}2,
\quad
\frac 52<q\le 3,
\label{rwttwo}
\end{align}
follow from
\eqref{L2Tlestimate} and
\eqref{Tlestimate}
 by a now standard  interpolation  argument due to Bourgain \cite{Bo1}
(see also the appendix in  \cite{CSWW}).
Of course \eqref{rwtone} and  \eqref{tildeT52estimate} imply
\eqref{rwt}. By a further interpolation (by the real method) we
can upgrade \eqref{rwttwo} to
\begin{equation}\label{strongtype}
\Big\|\sum_{2^l\leq \la^{1/3}}
T_{\la,l}\Big\|_{L^{p,q}(\bbR^2)\to L^{q}(\bbR^2)} \lc  \la ^{-2/q} ,\quad q=\frac{3p'}2,\quad
\frac 52<q< 3,
\end{equation} which
implies the analogous $L^p\to L^q$ inequality, and
we obtain \eqref{qgr52}, in the range $5/2<q< 3$. We note  that the case
$q\ge 3$
 (corresponding to  $p\le 2$) is already covered by the result in
\cite{GS1}.
Finally the inequality \eqref{ple52} follows by interpolation
between
 the $L^2(\bbR^2)$ bound $\|T_\la\|_{L^2\to L^2}=O(\la^{-5/6})$ and
the restricted weak type estimate \eqref{rwt}.

\medskip

\section{Bounds for  model cases}\label{models}
Consider the phase function defined in $\Bbb R^d$,
\begin{equation}\label{modelphi}\phi(x,y)=
\sum_{j=1}^{d-1}x_jy_j+\frac{(x_d-y_d)^3}{6}+ x_d\sum_{k=1}^{d-1}y_k^2,
\end{equation}
and let $\chi\in C_0^\infty(\bbR^d\times\bbR^d)$ be supported
near the origin.

We observe that $\rank \phi_{xy}''=d-1$, $\det \phi_{xy}''=
x_d-y_d$, and for $x_d=y_d$ the kernel of $d\pi_L$ is generated by
$\partial /\partial y_d$  and the kernel of  $d\pi_R$ is generated
by $\partial /\partial x_d$. Condition \eqref{piLfold} is
satisfied since $\phi_{x_dy_dy_d}=-1$ and condition
\eqref{piRfold} is satisfied since $\phi_{x_dy_dx_d}=1$. For each
$x$ the hypersurface $\cL_x=\{\phi_x'(x,y): \det\phi_{xy}=0\}$ is
just the paraboloid $\{(y',|y'|^2)\}$; thus condition
\eqref{curvcond} is satisfied.

Consider the operator $T_{\la,l}$ given by the localization  to
the set $\{|x_d-y_d|\approx 2^{-l}\}$. We now split $f=\sum f_m$
where $f_m(y)= \chi_{I_{m,l}}(y_d) f(y)$ and $I_{m,l}=[m2^{-l},
(m+1)2^{-l}]$. Then $T_{\lambda,l} f_m(x)$ vanishes if $x_d\notin
I_{m,l}^*:=I_{m-1,l}\cup I_{m,l}\cup I_{m+1,l}$. Thus
\begin{equation}
\label{fmtogether}
\|T_{\la,l}
f\|_q\leq C\Big(\sum_m\|T_{\la,l} f_m\|_q^q\Big)^{1/q}
\end{equation} and so it suffices to estimate
$T_{\la,l}f_m$.
Now we write
$$
T_{\la,l} f_m=\int_{I_{m,l}} e^{i\la(x_d-y_d)^3/6}T_{\la,l,y_d}
[f_m(\cdot,y_d)] dy_d,
$$ where for $g$ being defined on $\Bbb R^{d-1}$,
$$T_{\la,l,y_d} g(x',x_d)=
\int e^{i\la\Psi(x,y')}
 \chi(x,y',y_d)
g(y') dy'$$
and
\begin{equation}\Psi(x,y')= \sum_{j=1}^{d-1}x_jy_j+ x_d\sum_{k=1}^{d-1}y_k^2.
\label{psidef}
\end{equation}

Now let $d=2$. The phase function $\Psi$ is such that we can apply
 the Fefferman-Stein adjoint restriction theorem (\cite{Feff}),
or the more general Carleson-Sj\"olin theorem (\cite{CaSj}, \cite{H})
and obtain the estimates
$$
\|T_{\la,l,y_2} g\|_{L^q(\bbR^2)} \lc \lambda^{-2/q}
\|g\|_{L^p(\bbR) }, \quad q=3p', \ p<4;
$$
uniformly in $y_2\in I_{m,l}$. Then
\begin{align*}
\|T_{\la,l} f_m\|_{L^q(\bbR^2)} \lc &\int_{I_{m,l}}\|T_{\la,l,y_2}
[f_m(\cdot,y_2)] \|_{L^q(\bbR^2)} dy_2 \lc \lambda^{-2/q}
\int_{I_{m,l}} \|f_m(\cdot,y_2)\|_{L^p(\bbR)} dy_2
\\ &\lc 2^{-l/p'}
 \lambda^{-2/q}
\Big(\int_{I_{m,l}} \|f_m(\cdot,y_2) \|_{L^p(\bbR^2)}^p
dy_2\Big)^{1/p} \lc 2^{-l/p'}
 \lambda^{-2/q}
\|f_m\|_{L^p(\bbR^2)},
\end{align*}
and  \eqref{Tlpqestimate} is now implied by \eqref{fmtogether}.
The estimate \eqref{Tlpqestimatelim} follows in a similar way.
Moreover the bounds \eqref{Tlestimate} and \eqref{Tlestimatelim}
follow by using an endpoint $L^4$ bound of the Carleson-Sj\"olin
theorem.

\medskip
\noi{\it Higher dimensions.} A similar argument gives also a partial result
in higher dimensions.
Namely,  for the operator with model phase
 \eqref{modelphi} there is the bound
\begin{equation}\label{Tlataobound}
\big\|T_{\la}\big\|_{L^{p}(\bbR^d)\to L^{q}(\bbR^d)}\leq \la^{-d/q},
\quad  q\geq
\frac{d+1}d p', \quad q> \frac{2(d^2+d-1)}{d^2}.
\end{equation}
The range $q\ge 2(d+1)/d$ is covered by \cite{GS1}, and for
$2(d^2+d-1)/d^2<q<2(d+1)/d$ one can use
Tao's adjoint restriction estimate for paraboloids \cite{Tao}.
Indeed this estimate implies that the $L^r(\bbR^{d-1})\to L^s(\bbR^d)$
operator norm of
$T_{\la,l,y_d}$ is $O(\la^{-d/s})$, provided that
$s=\tfrac{d+1}{d-1} \, r'$ and  $ s>2(d+2)/d$.
By the above argument using H\"older's inequality in the $y_d$ variable
$$\|T_{\la,l} f\|_{s} \lc  2^{-l/r'} \lambda^{-d/s} \|f\|_r,
\quad s=\frac{d+1}{d-1} \, r', \quad \frac{2(d+2)}{d}<s< \frac {2(d+1)}{d-1}.
$$
We also use the $L^2\to L^2$ bound
\eqref{L2Tlestimate} and Bourgain's interpolation lemma. One deduces
 that the  operators
$\sum_{2^l\leq \la^{1/3}}T_{\la,l}$ map $L^{p,1}$ to
$L^{q,\infty}$ with norm $O( \la^{-d/q})$ if $ q\leq (d+1)p'/d$ and
$q> 2(d^2+d-1)/d^2$. By a further interpolation the
strong type $L^p\to L^q$ bound now follows in the same range; moreover there
 are similar bounds for $\widetilde T_\la$. Hence one obtains
\eqref{Tlataobound} in the full range.


We conjecture that this behavior remains true for general
oscillatory integral operators with folding canonical relations,
satisfying the elliptical curvature condition
\eqref{curvcond}. Well-known (hyperbolic) examples of Bourgain in
\cite{Bo2} may be adapted to show that an ellipticity condition is
in fact necessary here. We hope to pursue these  questions in a
subsequent paper.

It is conjectured that  the oscillatory integral operator $S_\la$
associated to the Carleson-Sj\"olin model phase $\Psi(x,y')$ as in
\eqref{psidef} has an $L^r(\bbR^{d-1})\to L^s(\bbR^d)$ operator
norm $O(\la^{-d/s})$ for $s\geq \frac{d+1}{d-1} r'$, $r<2d/(d-1)$.
The above analysis suggests that the bound \eqref{Tlataobound} for the
model case
might  be valid in the  range $q>(2d+1)/d$. Note that
$2(d^2+d-1)/d^2\geq (2d+1)/d$ for $d\geq 2$,
 with equality only  for $d=2$.

\medskip
\noi {\it One-sided fold conditions.}  Examples
suggest that the $L^p\to L^q$
estimates in Theorem \ref{main2d} for $p>5/2$ may hold merely
under the one-sided assumption \eqref{piLfold} and the curvature condition
\eqref{curvcond}.  This is in contrast to the $L^2$ estimates where
the bounds depend on finite type conditions on the projection $\pi_R$, see
\cite{GS1}, \cite{Com} and also the survey \cite{GSsurvey}.

A simple  example (where $\pi_R$ is maximally degenerate) is given by
\begin{equation}\label{modelphionesided}\psi(x,y)=
\sum_{j=1}^{d-1}x_jy_j+x_dy_d^2+ x_d\sum_{k=1}^{d-1}y_k^2.
\end{equation}
Now $\det \psi_{xy}=2y_d$ and a simple modification of the above
slicing argument shows that the $L^p\to L^q$ operator norm of the corresponding oscillatory integral operator is again $O(\la^{-d/q})$ if
$q\ge(d+1)p'/d$ and $q>2(d^2+d-1)/d^2$.
It would be interesting to know whether this result holds in general
under merely the  assumptions
\eqref{piLfold},
\eqref{curvcond}.  In \cite{GS1} it had been shown that this is the case in the range $q\ge 2(d+1)/d$.

\section{Estimation of $T_{\la,l}$ in two dimensions}
\label{proofpartone}
We shall now fix $l$ and various decompositions will depend on $l$
but this will not be indicated.
We shall estimate the square of $T_{\la,l}f $ and
bilinearize the problem as follows.
We split
$$ (T_{\la,l} f)^2= \sum_{m\geq 0} \fB^{m} (f,f)
$$
where
\begin{equation}\label{B0}\fB^{0} (f,\tf)(x)= \iint\chi_0(2^{-5+l/2}(y_1-\ty_1))
e^{i\la(\phi(x,y)+\phi(x,\ty))} \zeta_l(x,y)\zeta_l(x,\ty) f(y)
\tf(\ty) dy d\ty,
\end{equation}
and, for $m>0$,
\begin{equation}\label{Bm}
\fB^{m} (f,\tf)(x)= \iint\chi_1(2^{l/2-m-5}(y_1-\ty_1))
e^{i\la(\phi(x,y)+\phi(x,\ty))} \zeta_l(x,y)\zeta_l(x,\ty) f(y)
\tf(\ty) dy d\ty,
\end{equation}
and $$\zeta_l(x,y)=\chi(x,y)\chi_1(2^l(y_2-g(x,y_1))).$$
Notice that  the  sum in $m$  is extended over those  $m\ge 0$ with
$m<l/2 -C$ for large $C$ in view of the smallness of the support
of the cutoff function.

We shall show that

\begin{equation}
\label{Bpqestimate}
\|\fB^{0} (f,\tf)\|_{q/2}+\Big \|
\sum_{m>0}\fB^{m} (f,\tf)\Big\|_{q/2}
\lc 2^{-2l/p'}\la^{-4/q}
 \|f\|_p\|\tf\|_p \quad\text{ for $1\leq p< 4$, $q=3p'$};
\end{equation}
moreover
\begin{equation}
\label{B0estimate}
\|\fB^{0} (f,\tf)\|_2 \lc 2^{-3l/2}\la^{-1} (\log \la )^{1/2}
 \|f\|_4\|\tf\|_4
\end{equation}
and
\begin{equation}
\label{B0mestimate}
\Big\|\sum_{m>0}\fB^{m} (f,\tf)\Big\|_2 \lc 2^{-3l/2}(1+l)^{1/2}\la^{-1}
 \|f\|_4\|\tf\|_4.
\end{equation}

In what follows we shall  estimate the expression
$\sum_{m>0}\fB^m(f,f)$ for $m>0$ and give the modifications for
$\fB^0(f,f)$ in \S\ref{B0bound}.

The principal objective of our approach is to reduce matters
to  an $L^2$  estimate for some well-localized  operators
(termed $\fS\equiv \fS^{m\mu b\tb}_{a\nu}$ below), for which
one can use arguments for  model cases considered
in \S\ref{models}.
The idea is to estimate such a localized  operator
$\fS$, by freezing the variables $y_2, \widetilde y_2$,
and to take advantage of the small support by using H\"older's inequality.
It was possible to  implement this idea ``globally'' in  ``rigid''
model cases, such as  $g(x,y_1)= x_2+ x_1-y_1$, but
this global approach does not seem to  work for  general  $g$.

It seems natural to decompose for fixed $x$
 the set $\{y:|y_2-g(x,y_1)|\approx 2^{-l}\}$ into rectangular pieces of
size $\approx 2^{-l/2}\times 2^{-l}$. In order to be useful this
decomposition should be stable under perturbations in $x$
but since  $g(x,y_1)$ varies in $x$ we need a decomposition in $x$ as well.
For a situation in which we can use the idea of freezing $y_2$
we may consider the case that
 $y_1,\tilde y_1$ are supported in intervals  $I,  \widetilde I$ of length
$2^{-l/2}$ and $x$  is localized to certain rectangles $R$ of size
$2^{-l/2}\times 2^{-l}$ so that for  $(x,y_1,\widetilde y_1)\in
R\times I\times \widetilde I$ the expressions $g(x,y_1)$ and
$g(x,\tilde y_1)$ vary by no more than $2^{-l}$. This works well
if the distance of $I$ and $\widetilde I$ is not much more than
their length, namely $2^{-l/2}$. The  rectangles in $x$-space are
not supposed to change orientation while their centers vary over a
cube of sidelength $2^{-l/2}$ and the geometry of $\cL$ suggests
that the long sides become perpendicular to $\nabla_x g$ at the
centers of the cubes.

When estimating the $\fB^m(f,\widetilde f )$ we are in the
situation where $|y_1-\widetilde y_1|\approx 2^{-l/2+m+5}$ for $0<
m<l/2$. We would then
 like to make a similar decomposition of
 intervals $I\times \widetilde I$ in $(y_1,\widetilde y_1)$ space
and rectangles $R$ in $x$ space. The requirement that {\it both}
 expressions $g(x,y_1)$ and $g(x,\tilde y_1)$
vary by no more than $2^{-l}$ in $R\times I\times\widetilde I$
is  now harder  to satisfy and we need to choose a finer decomposition,
namely we choose the intervals $I$, $\widetilde I$ to be of length
$2^{-l/2-m}$ ({\it cf.} \eqref{ommba} below) and make a decomposition
in terms in $x$ into cubes  of the sidelength
$2^{-l/2-m}$  ({\it cf.}
\eqref{uma} below). This is somewhat reminiscent of a
situation in \cite{Se}.
Moreover we  decompose each cube in $x$ space
into smaller rectangles of sidelengths
$2^{-l/2-m}$ and $2^{-l}$ and
the longer sides  are  perpendicular to  $\nabla_x g$
 at the centers of the cubes
({\it cf.} \eqref{ummuanu} below). The geometry is now such that
the orientation of the rectangles is essentially the same when
$y_1$ varies over $I$,
$\widetilde y_1$ varies over $\widetilde I$ and $x$ varies over $R$.
Note that the rectangles become essentially cubes of length $2^{-l}$
 when the distance  of the intervals  is $\approx 1$;
 in this case
the length of the intervals $I$, $\widetilde I$ is also $\approx 2^{-l}$.

We now formally define these  decompositions  in the $x$ and $y$
variables. Let $\eta$ be a $C^\infty$ function supported in
$(-3/4,3/4)$ so that
$\sum_{n\in\bbZ} \eta^2(s-n)\equiv 1.$
For
each $b\in
2^{-l/2-m}\mathbb{Z}$ let $J_b^m$ be the interval of
length $2^{-l/2-m+1}$ centered at $b$. We let $\cP^m$ be the set
of all pairs $(b,\tb)$ with the property that $J_b^m\times
J_\tb^m$ intersects the support of the cutoff function
$(y_1,\ty_1)\mapsto \chi_1(2^{l/2-m-5}|y_1-\ty_1|)$, so that for
$(b,\tb)\in \cP^m$ the numbers $b$, $\tb$ are
$C2^{-l/2+m}$-separated. We may split $\cP^m=\cup_{\mu} \cP^m_\mu$
where $\mu$ ranges over $2^{-l/2+m}\bbZ$, the families $\cP^m_\mu$
are disjoint and of cardinality $O(2^{4m})$ and we have
$$|b-\mu|\leq 2^{-l/2+m+1}, \quad|\tb-\mu|\leq 2^{-l/2+m+1}, \quad\text {for } (b,\tb)\in \cP^m_\mu.$$

We also decompose in $x$ space, using two parameters $a,\nu$. The
parameter $a$ will range over points in $2^{-l/2-m}\bbZ^2$. For
$\mu\in 2^{-l/2+m}\bbZ$ and   $\nu\in \bbZ$ (typically $|\nu|\leq
C2^{l/2-m}$)  we set
$$\tau^\mu_\nu a:=a+2^{-l}\nu \tfrac{\nabla_xg(a,\mu)}{|\nabla_xg(a,\mu)|}
$$ and
\begin{align}
u^{m}_{a}(x)
\label{uma}&:= \eta(2^{l/2+m}(x-a))
\\
u^{m\mu}_{a\nu}(x) \label{ummuanu} &:= \eta(2^{l/2+m}(x-a))
\eta(2^{l}(\inn{\tfrac{\nabla_xg(a,\mu)}{|\nabla_xg(a,\mu)|}}{x-\tau^\mu_\nu
a})).
\end{align}
Moreover for $(b,\tb)\in
\cP^m_\mu$,
we make the following definitions. Let $\beta_0\in C^\infty_0(\bbR) $ be equal to $1$ near $0$ so that
 $\chi_1\beta_0=\chi_1$. Set
\begin{align}
\label{ommba}
\omega^{mb}_{a}(y)
&=\eta(2^{l/2+m}(y_1-b))
\beta_0(C^{-1}2^{l/2+m}(y_2-g(a,y_1)))
\\
\label{ommmubanu}
\omega^{m\mu b}_{a\nu}(y)
&:=\omega^{mb}_{a}(y)
\beta_0(C^{-1}2^{l}(y_2-g(\tau_\nu^\mu a,y_1)))
\end{align}
and finally
\begin{equation} \label{Ommbtb}
\Omega^{mb\tb}_{a}(y,\ty): =\omega^{mb}_{a}(y) \omega^{m\tb}_{a}(\ty),
\qquad
\Omega^{m\mu b\tb}_{a\nu}(y,\ty)
 :=\omega^{m\mu b}_{a\nu}(y)
 \omega^{m\mu\tb}_{a\nu}(\ty).
\end{equation}

%

For locally integrable  functions $F$ defined on $\bbR^4$ and
$m>0$ we set
\begin{multline}\fS^{m\mu b\tb}_{a\nu}F(x)=
(u^{m\mu}_{a\nu}(x)) ^2
\\ \times
\iint
\chi_1(2^{l/2-m-5}(y_1-\ty_1))
\zeta_l(x,y)  \zeta_l(x,\ty)
e^{i\la(\phi(x,y)+\phi(x,\ty))}
\Omega^{m\mu b\tb}_{a\nu}(y,\ty)
F(y,\ty)
dy d\ty;
\end{multline}
and for $m=0$ we use a similar  definition with the modification that
$\chi_1(2^{l/2-m-5}(y_1-\ty_1))$ is replaced  by
$\chi_0(2^{l/2-5}(y_1-\ty_1))$.

Typically the operator $\fS^{m\mu b\tb}_{a\nu} $ should be acting
on the   function
\begin{equation}\label{tensorF}
(y,\ty)\mapsto f^{m\mu b}_{a\nu}\otimes \tf^{m\mu \tb}_{a\nu}
(y,\ty):= \omega^{m\mu b}_{a\nu}(y)f(y) \omega^{m\mu
\tb}_{a\nu}(\ty)\tf(\ty)
\end{equation}
when $(b,\tb)\in\cP^m_\mu$. Indeed,
in view of the condition
$\sum_{n\in\bbZ} \eta^2(s-n)\equiv 1$, we have
$$\fB^m(f,f)=\sum_{a,\nu} \sum_{\mu}
\sum_{(b,\tb)\in \cP^m_\mu} \fS^{m\mu b\tb}_{a\nu} ( f^{m\mu
b}_{a\nu}\otimes \tf^{m\mu \tb}_{a\nu}).$$

\medskip

In \S\ref{proofoffirstprop} we shall prove the following
inequality concerning  vector-valued functions, which combines
various orthogonality arguments with the individual estimates for
the operators $\fS^{m\mu b\tb}_{a\nu}$. As we shall see, the proof
relies on ideas related to the  Carleson-Sj\"olin  theorem.

\begin{proposition}\label{firstprop}
For $2\leq r\leq \infty$
\begin{multline}\label{firstpropeq}
\Big\|\sum_{0<m<l/2}\sum_{\mu\in 2^{-l/2+m}\bbZ} \sum_{(b,\tb)\in \cP_\mu^m} \sum_{a,\nu}
\fS^{m\mu b\tb}_{a\nu}
F^{m\mu b\tb}_{a\nu}\Big\|_r \lc\\ \la^{-2/r}
\Big(
\sum_{0<m<l/2}2^{-(m+3l/2)\frac{r'}{r} }
\sum_\mu \sum_{(b,\tb)\in \cP_\mu^m}
\Big(
\sum_{a,\nu} \Big(\iint|F^{m\mu b\tb}_{a\nu}(y,\ty)|^{r'}dyd\ty\Big)^{\frac{r}{r'}}
\Big)^{\frac{r'}{r} }\Big)^{\frac 1{r'}}.
\end{multline}
\end{proposition}

We shall have to choose the functions $F^{m\mu b\tb}_{a\nu}$
carefully in order to take full advantage of Proposition
\ref{firstprop}. As mentioned above we would like to let
$\fS^{m\mu b\tb}_{a\nu}$ act on
 the function \eqref{tensorF}.
However we shall have to exploit finer  frequency localization
properties  of the operator $\fS^{m\mu b\tb}_{a\nu}$. Split
\begin{equation}
f^{m\mu b}_{a\nu}= \cL f^{m\mu b}_{a\nu}(y)+
\cE f^{m\mu b}_{a\nu}(y)
\end{equation}
where
\begin{equation}
\label{Lfrequencyloc}
\cL f^{m\mu b}_{a\nu}(y)=
\omega^{m\mu b}_{a\nu}(y)
\int f(z_1,y_2) \omega^{mb}_{a}(z_1,y_2)
\int \beta_2\big(\tfrac{\la\phi_{y_1}(a,b, g(a,b)) +\eta}
{\la 2^{-l/2-m}}\big)
e^{i\eta(y_1-z_1)} \frac{d\eta}{2\pi} dz_1
\end{equation}
and
\begin{equation}
\label{Efrequencyloc}
\cE f^{m\mu b}_{a\nu}(y)=
\omega^{m\mu b}_{a\nu}(y)\int f(z_1,y_2)
\omega^{mb}_{a}(z_1,y_2)
\int \Big(1-\beta_2\big(\tfrac{\la\phi_{y_1}(a,b, g(a,b)) +\eta}
{\la 2^{-l/2-m}}\big)\Big)
e^{i\eta(y_1-z_1)} \frac{d\eta}{2\pi} dz_1;
\end{equation}
here
the function $\beta_2$ is supported in
the union of $(-2C_0, -(2C_0)^{-1})$ and
$((2C_0)^{-1}, 2C_0)$ and $\beta_2(s)=1$ if $|s|\in [C_0^{-1}, C_0]$,
for suitably large $C_0$,
and the integral in the definition
of $\cE f^{m\tb}_{a\nu}$ is to be interpreted as an oscillatory integral.
Now
$$\fS^{m \mu b\tb}_{a\nu}( f^{m\mu b}_{a\nu}\otimes \tf^{m\mu \tb}_{a\nu})=
\fS^{m\mu b\tb}_{a\nu}\big( (\cL f^{m\mu b}_{a\nu} + \cE f^{m\mu
b}_{a\nu}) \otimes (\cL\tf^{m\mu \tb}_{a\nu}+\cE\tf^{m\mu
\tb}_{a\nu})\big),
$$
but only the contribution of
$\cL f^{m\mu b}_{a\nu}\otimes
\cL \tf^{m\mu \tb}_{a\nu}$ is relevant:

\begin{lemma} \label{yerrorterms}
There are the pointwise bounds
\begin{align}
&\fS^{m\mu b\tb}_{a\nu}\big( \cE f^{m\mu b}_{a\nu}
\otimes \cL\tf^{m\mu \tb}_{a\nu}\big)= O((\la 2^{-2l})^{-N})
\\
&\fS^{m\mu b\tb}_{a\nu}\big( \cL f^{m\mu b}_{a\nu}
\otimes \cE\tf^{m\mu \tb}_{a\nu}\big)
= O((\la 2^{-2l})^{-N})
\\
&\fS^{m\mu b\tb}_{a\nu}\big( \cE f^{m\mu b}_{a\nu}
\otimes \cE\tf^{m\mu \tb}_{a\nu}\big)
= O((\la 2^{-2l})^{-N})
\end{align}
\end{lemma}
\begin{proof}[Sketch of proof.]

We only consider the term $\fS^{m\mu b\tb}_{a\nu}\big( \cE f^{m\mu
b}_{a\nu} \otimes \cL\tf^{m\mu \tb}_{a\nu}\big)$; the others are
handled similarly.

The analysis leads to the estimation of  oscillatory integrals of the form
$$
\iint K(x,z_1,y,\tilde z_1, y) e^{i (\la \phi(x,y_1,y_2)+\eta
y_1)}\Big(1-\beta_2\big(\tfrac{\la\phi_{y_1}'(a,b, g(a,b)) +\eta}
{\la 2^{-l/2-m}}\big)\Big) dy_1 d\eta
$$
where $K$ is a function satisfying $\partial_{y_1}^\alpha K=O(
2^{\alpha l})$ which vanishes for $|x-a|\gc 2^{-m-l/2}$,
$|y_1-b_1|\gc 2^{-m-l/2}$, $|y_2-g(x,y_1)|\gc 2^{-l}$. The $y_1$
derivative of the phase is then
$$
\la\phi_{y_1}(x,y_1,y_2)+\eta= \la\phi_{y_1}'(a,b, g(a,b)) +\eta +
O(\la 2^{-l/2-m}).
$$
As we assume the constant $C_0$ in the definition of the
cutoff function $\beta_2$  to be large we see that
$$\big|\la\phi_{y_1}(x,y_1,y_2)+\eta\big|\gc
|\la\phi_{y_1}'(a,b, g(a,b)) +\eta |\gg \la 2^{-l/2-m}.
$$
The assertion then follows by an integration by parts with respect to the $y_1$ variable.
\end{proof}

We shall need an orthogonality property of the $\cL f^{m\mu
b}_{a\nu}$. Set
\begin{equation} \label{fmb}
f^{mb}(y)=\chi_{[b-2^{-l/2-m},b+2^{-l/2-m}]}(y_1) f(y).
\end{equation}
\begin{lemma}\label{lpLplemma} For $ p\geq 2$ and fixed $m, b, \mu$
\begin{equation}
\label{lpLpinequality}
\Big(\sum_{a,\nu}\|\cL f^{m\mu b}_{a\nu}\|_p^p\Big)^{1/p}
\lc \|f^{mb}\|_p
\end{equation}
uniformly in $m,\mu, b$.
\end{lemma}

\begin{proof}

First note that for fixed $m,\mu,b, a$
the supports of the functions
$$\beta_0(C^{-1}2^{l/2}(y_2-g(\tau_\nu^\mu a,y_1)))$$
have uniformly  bounded
 overlap on the  support of $\omega^{m b}_{a}$.

Define
\begin{equation}L^{m b}_{a} f(y)= \omega^{mb}_a(y)
\int f(z_1,y_2) \omega^{mb}_{a}(z_1,y_2)
 \frac{1}{2\pi}
\int \beta_2\big(\tfrac{\la\phi_{y_1}'(a,b, g(a,b)) +\eta} {\la
2^{-l/2-m}}\big) e^{i\eta(y_1-z_1)} d\eta\,
dz_1.
\end{equation}

Then the left hand side of
\eqref{lpLpinequality} is dominated by a constant times the left hand side of the following inequality
\begin{equation}
\label{newlpLpinequality}
\Big(\sum_{a}\|L^{m b}_a f^{mb}\|_p^p\Big)^{1/p}
\lc \|f^{mb}\|_p
\end{equation}
which we now prove. It is easy to see that the operators
$L^{mb}_a$ are uniformly bounded on $L^\infty$ and an
interpolation argument reduces the proof of
\eqref{newlpLpinequality} to the case $p=2$.

In order to complete this proof it suffices to check that
\begin{equation}
\label{orthogonality}
\big\|
L^{m b}_a (L^{m b}_{a'})^*
\big\|_{L^2\to L^2}
\lc 2^{-|n-n'|} \quad \text{ if } a= 2^{-m-l/2} n, \
a'= 2^{-m-l/2} n', \ |n-n'|\geq M,
\end{equation}
for suitably large $M$.

The kernel of $L^{m b}_a (L^{m b}_{a'})^*$ is given by
\begin{multline} \label{schwartzkernel}
\cK(y,y')=\delta(y_2-y_2')
\omega_a^{mb}(y) \omega_a^{mb}(y')
\int \omega_a^{mb}(z_1,y_2) \omega_a^{mb}(z_1, y_2')
\\ \times\Big\{
\iint e^{i((y_1-z_1)\eta-(y_1'-z_1)\eta')}
\beta_2\big(\tfrac{\la\phi_{y_1}(a,b, g(a,b)) +\eta} {\la
2^{-l/2-m}}\big) \beta_2\big(\tfrac{\la\phi_{y_1}(a',b, g(a',b))
+\eta'} { \la 2^{-l/2-m}}\big) \frac{d\eta}{2\pi}
\frac{d\eta'}{2\pi}\Big\}\,dz_1
\end{multline}

We shall now fix $a'$ and solve the equation
$g(x_1,x_2, b)= g(a_1',a_2',b)$ in $x_2$; this can be done by the
implicit function theorem since
$$g_{x_2}(x,y_1)=-\frac
{\phi_{x_2y_2x_2}(x,y_1,g(x,y_1))}
{\phi_{x_2y_2y_2}(x,y_1,g(x,y_1))} +o(1)
$$
where $o(1)$ is a quantity which by \eqref{phixdy}
vanishes at the  reference point
$P$ and the two sided fold assumption
(\eqref{phixdydyd}, \eqref{phixdxdyd})  implies $g_{x_2}\neq 0$.
For later reference we also note that
\begin{equation}
\label{smallnessofgx1}
 g_{x_1}(x,y_1)=-\frac{\phi_{x_2y_2x_1}(x,y_1,g(x,y_1))}
{\phi_{x_2y_2y_2}(x,y_1,g(x,y_1))} +o(1) =o(1)
\end{equation}
which follows from \eqref{phixdy} and \eqref{mixedxdterm}.

Let thus $h(x_1,a',b)$ denote the unique solution satisfying
\begin{equation}
\begin{aligned}
g(x_1, h(x_1,a',b),b)
&=g(a',b),
\\
h(a_1',a',b)&=a_2'.
\end{aligned}
\end{equation}
Then
\begin{equation}\begin{aligned}
g(a,b)-g(a',b)&=g(a,b)- g(a_1, h(a_1,a',b),b)+
g(a_1, h(a_1,a',b),b)-g(a',b)
\\
&=\cC(a,a',b) (a_2-h(a_1,a',b))+ O(a_1-a_1')
\end{aligned}\end{equation}
where $\cC(a,a',b)\neq 0$.
Thus if  for some small constant $c_0$ and some large constant $C_0$
\begin{equation}\label{caseone}
|a_1-a_1'|\leq c_0|a_2-h(a_1,a',b)| \quad \text{ and }
|a_2-h(a_1,a',b)| \geq C_0 2^{-m-l/2},
\end{equation}
then
$|g(a,b)-g(a',b)|\geq C_1 2^{-m-l/2}$ for still large $C_1$ and therefore
we have
\begin{equation} L^{m b}_a (L^{m b}_{a'})^*\,=\,0\end{equation}
in the case \eqref{caseone}.

In the relevant opposite case we assume that
\begin{equation}\label{casetwo}
|a_1-a_1'|\geq c_0|a_2-h(a_1,a',b)| \quad \text{ and }
|a_1-a_1'|\geq C_1 2^{-m-l/2}
\end{equation}
Since by \eqref{smallnessofgx1} we also have $h_{x_1}=o(1)$ it
follows that in the present case \eqref{casetwo} we have
$|a_1-a_1'|\gg |a_2-a_2'|$ and therefore we can estimate  with
$P_{a'b}=(a',b,g(a',b))$
\begin{equation*}
\begin{aligned}
&\phi_{y_1}(a, b, g(a,b))-
\phi_{y_1}(a', b,g(a',b))
\\&=
\phi_{y_1}(a,b, g(a,b))-\phi_{y_1}(a_1, h(a_1,a',b),
b,g(a_1,h(a_1,a',b),b))
\\&\quad+ \phi_{y_1}(a_1, h(a_1,a',b), b,g(a_1,h(a_1,a',b),b))-
\phi_{y_1}(a',b, g(a',b))
\\
&=\phi_{x_1y_1}(P_{a'b}) (a_1-a_1')
+O(\eps|a_1-a_1'|)+
O(\eps
|a_2-h(a_1,a',b)|)
\end{aligned}
\end{equation*}
if the support of the initial cutoff function has diameter $\leq \eps$.
Here, in order to get the $\eps$-bound,  we used the assumption \eqref{phixdy}.
Thus in case \eqref{casetwo} we get
\begin{equation*}
|\phi_{y_1}(a, b, g(a,b))- \phi_{y_1}(a', b,g(a',b))|\geq
C|a_1-a_1'|\approx |a-a'|.
\end{equation*}
Hence by an integration by parts in the $z_1$ variable we gain
negative powers of $$2^{-l}|\eta-\eta'|\gc 2^{-m-3l/2}\la|n-n'|$$
if $|n-n'|$ is large; this
 is more than enough to prove the required almost orthogonality property.
\end{proof}

Now applying Proposition \ref{firstprop}  to the functions
\begin{equation}F^{m\mu b\tb}_{a\nu}(y,\ty)=
 \cL f^{m\mu b}_{a\nu} (y)  \cL \tf^{m\mu \tb}_{a\nu} (\ty)
\end{equation}
we reduced  matters to

\begin{proposition} \label{secondprop}
For $2<r<\infty$, $r=3p'/2$,
\begin{multline}\label{secondpropestimate}
\Big(
\sum_{0<m<l/2}\sum_\mu \sum_{(b,\tb)\in \cP^{m}_{\mu}  }
2^{-(m+3l/2)\frac{r'}{r} }\Big(
\sum_{a,\nu} \Big(\iint|\cL f^{m\mu b}_{a\nu} (y)
\cL \tf^{m\mu \tb}_{a\nu} (\ty)
|^{r'}dyd\ty\Big)^{\frac{r}{r'}}
\Big)^{\frac{r'}{r} }\Big)^{\frac 1{r'}}
\\
 \lc 2^{-2l/p'}\|f\|_p\|\tf\|_p.
\end{multline}
For $r=2$, $p=4$ the  left hand side of
\eqref{secondpropestimate}
is dominated by
$$C (1+l)^{1/2} 2^{-3l/2}\|f\|_4\|\tf\|_4 .$$
\end{proposition}

\begin{proof}

Each  $\cL f^{m\mu b}_{a\nu}$ is supported on a parallelogram of
area $\lc 2^{-m-3l/2}$ and so by H\"older's inequality we can
 estimate for $r\geq 2$ the left hand side of
\eqref{secondpropestimate}  by
\begin{align}
&\Big(
\sum_{m<l/2} \sum_\mu\sum_{(b,\tb)\in \cP_\mu^m}  2^{-(m+3l/2)\frac{r'}{r} }
2^{-(m+3l/2)(2-\frac {2r'}p)}
\Big(
\sum_{a,\nu}
\|\cL f^{m\mu b}_{a\nu}\|_p^{r}
\|\cL \tf^{m\mu \tb}_{a\nu}\|_p^{r}\Big)^{\frac{r'}{r}}
\Big)^{\frac{1}{r'}}
\notag
\\
&\label{expr1} \lc\Big( \sum_{m<l/2}
2^{-(m+3l/2)(\frac{r'}{r}+2-\frac {2r'}p)}\sum_{(b,\tb)\in \cP^m}
\|f^{mb}\|_p^{r'} \|\tf^{m\tb}\|_p^{r'}\Big)^{\frac 1{r'}}
\end{align}
where $f^{mb}$ is an in \eqref{fmb}.
For the last inequality we have used the Cauchy-Schwarz inequality, the
fact that $2r\ge p$ (which follows from our assumptions on $p$ and $r$), the
 embedding $\ell^p\subset \ell^{2r}$ and \eqref{lpLpinequality} of Lemma \ref{lpLplemma}.

Now
 let
$I_n=[n2^{-l/2}, (n+1)2^{-l/2})$  for $n\in \bbZ$; then $I_n$ contains $\approx 2^{m}$ numbers $b\in 2^{-m-l/2}\bbZ$ and we  dominate
\eqref{expr1} by a constant times
\begin{align}
&\Big(
\sum_{m<l/2} \sum_{(n,\tn)\in \bbZ^2 \atop{ |n-\tn|\approx 2^{m}}}
2^{-(m+3l/2)(\frac{r'}{r}+2-\frac {2r'}p)}2^{2m(1-\frac {r'}p)}
\\&\quad\quad\quad\quad\times
\Big(
\sum_{b\in I_n\cap\atop{ 2^{-m-l/2}\bbZ} }
\|f^{m b}\|_p^{p}\Big)^{\frac {r'}p}
\Big(
\sum_{\tb\in I_\tn\cap\atop{ 2^{-m-l/2}\bbZ} }
\|\tf^{m \tb}\|_p^{p}\Big)^{\frac {r'}p}
\Big)^{\frac{1}{r'}}
\notag
\\
&\lc
2^{-2l/p'}
\Big(\sum_m
\sum_{(n,\tn):\atop
{|n-\tn|\approx 2^{m}}}|n-\tn|^{-r'+1}
\|f\|_{L^p(I_n\times\bbR)}^{r'}
\|\tf\|_{L^p(I_\tn\times\bbR)}^{r'}
\Big)^{\frac{1}{r'}};
\label{expr2}
\end{align}
here we have used that $1/r+2/r'-2/p= 4/(3p')$ in view of the assumption
$r=q/2=3p'/2$.

Let $\beta\in (0,1)$ and define  for a sequence $\fa$  the discrete analogue of the standard  fractional integral
$$[I^\beta\fa]_n=\sum_{\tn}|n-\tn|^{\beta-1}\fa_\tn.$$ Now the condition
$r=3p'/2$, is equivalent with
$2-r'=\frac{1}{p/r'}-\frac{1}{(p/r')'}$ so that for
$2<r<\infty$, $r=3p'/2$ the operator
$I^{2-r'}$ maps $\ell^{p/r'}\to \ell^{(p/r')'}$. We apply this with
$\fa_n=\|f\|_{L^p(I_n\times \bbR)}^{r'}$
and also set
$\tilde \fa_n=\|\tf\|_{L^p(I_n\times \bbR)}^{r'}$. Then
the expression \eqref{expr2} is bounded by
\begin{equation*}
C 2^{-2l/p'}\Big(\sum_n \tilde \fa_n [I^{2-r'}\fa]_n\Big)^{1/r'}
\end{equation*} and we argue as in H\"ormander \cite{H} to get
\begin{equation*}
\Big(\sum_n \big|\tilde \fa_n [I^{2-r'}\fa]_n\big|\Big)^{1/r'}\leq
 \|\tilde \fa\|_{p/r'}^{1/r'}\|I^{2-r'}\fa\|_{(p/r')'}^{1/r'}
\lc  \|\fa\|_{p/r'}^{1/r'}
\|\tilde \fa\|_{p/r'}^{1/r'}
\lc \|f\|_p\|\tf\|_p.
\end{equation*}

The case $r=2$, $p=4$ is similar, except that the expression \eqref{expr2}
is now estimated using a simple convolution inequality for each fixed $m$ and the sum over $m$ introduces the logarithmic term.
\end{proof}

\section{Proof of Proposition \ref{firstprop}}
\label{proofoffirstprop}

We prove
inequality \eqref{firstpropeq} by interpolation between the extreme cases $r=2$ and $r=\infty$.
The case $r=\infty$
 is
\begin{multline*}
\Big\|\sum_{0<m<l/2}\sum_{\mu\in 2^{-l/2+m}\bbZ} \sum_{b,\tb\in \cP_\mu^m} \sum_{a,\nu}
\fS^{m\mu b\tb}_{a\nu}F^{m\mu b\tb}_{a\nu}\Big\|_\infty \lc
\\
\sum_{0<m<l/2}
\sum_\mu \sum_{(b,\tb)\in \cP_\mu^m}
\sup_{a,\nu} \iint|F^{m\mu b\tb}_{a\nu}(y,\ty)|dyd\ty.
\end{multline*}
This is immediate; one uses for fixed $m,\mu$
the almost disjointness of the cutoff functions
$u^{m\mu}_{a\nu}$ in \eqref{ummuanu}.

For the remainder of this section we consider the case $r=2$ which
is
\begin{multline} \label{caser=2}
\Big\|\sum_{0<m<l/2}\sum_{\mu\in 2^{-l/2+m}\bbZ}
 \sum_{(b,\tb)\in \cP_\mu^m} \sum_{a,\nu}
\fS^{m\mu b\tb}_{a\nu} F^{m\mu b\tb}_{a\nu}\Big\|_2 \lc\\
\frac{1}{\la} \Big( \sum_{0<m<l/2}2^{-(m+3l/2)} \sum_\mu
\sum_{(b,\tb)\in \cP_\mu^m} \sum_{a,\nu} \iint|F^{m\mu
b\tb}_{a\nu}(y,\ty)|^{2}dyd\ty \Big)^{\frac 1{2}}.
\end{multline}

%
\subsection{The four steps in the proof}
We need to use  various orthogonality lemmata.

\begin{lemmasub}\label{firstortho}
For each $N\in\mathbb{N}$
\begin{align}\label{m-orthogonality}
\Big\|\sum_{0<m<l/2}&\sum_{\mu\in 2^{-l/2+m}\bbZ} \sum_{(b,\tb)\in
\cP_\mu^m} \sum_{a,\nu} \fS^{m\mu b\tb}_{a\nu}F^{m\mu
b\tb}_{a\nu}\Big\|_2 \lc
\\
&\Big(\sum_{m<l/2}
\sum_{a}
\Big\|\sum_{\mu\in 2^{-l/2+m}\bbZ\atop{
(b,\tb)\in \cP_\mu^m}} \sum_\nu
\fS^{m\mu b\tb}_{a\nu}F^{m\mu b\tb}_{a\nu}\Big\|_2^2\Big)^{\frac 12}
\notag
\\&+2^{-5l/4}(\la2^{-2l})^{-N}
\Big(\sum_{m,\mu,b,\tb,a,\nu}2^{-m(2N-1)}\|F^{m\mu
b\tb}_{a\nu}\|_2^2 \Big)^{\frac 12}. \notag
\end{align}
\end{lemmasub}

\begin{lemmasub}\label{1.5ortho}
For each $N\in\mathbb{N}$
\begin{multline}\label{mu-orthogonality}
\Big\|\sum_{\mu\in 2^{-l/2+m}\bbZ} \sum_{(b,\tb)\in \cP_\mu^m} \sum_{\nu}
\fS^{m\mu b\tb}_{a\nu}F^{m\mu b\tb}_{a\nu}\Big\|_2 \lc\\
\Big(\sum_{\mu\in 2^{-l/2+m}\bbZ}  \sum_{\nu} \Big\|
\sum_{(b,\tb)\in \cP_\mu^m} \fS^{m\mu b\tb}_{a\nu}F^{m\mu
b\tb}_{a\nu}\Big\|_2^2\Big)^{1/2}
+2^{-3l/2}(\la2^{-2l+m})^{-N}\Big( \sum_{\mu,b,\tb,\nu}\|F^{m\mu
b\tb}_{a\nu}\|_2^2\Big)^{1/2},
\end{multline}
uniformly in $m$ and $a$.
\end{lemmasub}

\begin{lemmasub}\label{secondortho}
For each $N\in\mathbb{N}$
\begin{multline}\label{b-orthogonality}
\Big\| \sum_{(b,\tb)\in \cP_\mu^m}
\fS^{m\mu b\tb}_{a\nu}F^{m\mu b\tb}_{a\nu}\Big\|_2 \lc\\
\Big(\sum_{(b,\tb)\in \cP_\mu^m} \big\|
\fS^{m\mu b\tb}_{a\nu}F^{m\mu b\tb}_{a\nu}\big\|_2^2\Big)^{1/2}
+2^{-3(2m+3l)/4}
 (\la2^{-m-3l/2})^{-N}
\Big(\sum_{b,\tb}\|F^{m\mu b\tb}_{a\nu}\|_2^2\Big)^{1/2},
\end{multline}
uniformly in $m>0,\mu, a$ and $\nu$.
\end{lemmasub}
In view of our condition  $2^l\le \la^{1/3}$ the precise error
bounds in the above lemmata will be unimportant.

These three estimates reduce matters to the uniform $L^2$ bounds for
the operators $\fS^{m\mu b\tb}_{a\nu}$:
%

\begin{propositionsub}\label{fixedfSlemma}
The estimate
\begin{equation}\label{fixedfSest}
\big\|
\fS^{m\mu b\tb}_{a\nu}
F\big\|_2 \lc 2^{-3l/4-m/2} \la^{-1}\|F\|_2,
\end{equation}
holds with bounds uniform in $m>0,\mu, a,\nu$ and $(b,\tb)\in
\cP^m_\mu$.
\end{propositionsub}

Inequality \eqref{firstpropeq} for $r=2$  is  an immediate consequence of
Lemmata  \ref{firstortho},  \ref{1.5ortho},  \ref{secondortho}
 and Proposition \ref{fixedfSlemma}; we take into account that $2^{2m}\leq 2^l$ and $2^l\leq \la^{1/3}$.

\subsection{Preliminary considerations}
We first state  a more or less standard result on oscillatory integrals,
for which we include a sketch of the proof for completeness.
\begin{lemmasub}
\label{basicestimate}
 Let $(x,y)\mapsto \Psi(x,y)$ be a smooth real valued phase
 function, defined in a domain $D\subset \Bbb R^d\times \Bbb R^d$
 so that $\rank (\Psi_{x'y'})= d-1$ in $D$  and so that
we have uniform bounds for the derivatives  of $\Psi$ in $D$;
\textit{i.e. },
\begin{equation}\label{unifupper}|\partial_{x,y}^\alpha\Psi|\leq C_\alpha,
\end{equation} for all $|\alpha|\leq 4d$.
Let $\lambda\gg 1$ and $\delta\geq \lambda^{-1/3}$. Let
$P^o=(x^o,y^o)$ and $Q_\delta(P^o)=\{(x,y):|x-x^o|\leq \delta,
|y-y^o|\leq \delta.\}$. Suppose that for some $C_1>0$
\begin{equation}\label{uniflower}C_1^{-1}\delta\leq |\det \Psi_{xy}|\leq C_1 \delta \text{ for } (x,y)\in
Q_\delta(P).\end{equation} Let $a$ be supported in $Q_\delta(P^o)$
and assume that
\begin{equation} \label{C-V-symbol} |\partial_{x,y}^\alpha a|\leq C_\alpha
(\lambda \delta)^{|\alpha|/2}
\end{equation}
for all multiindices $\alpha$. Define the operator $\cJ_\lambda$ by
$$\cJ_\la f(x)=\int e^{i\la \Psi(x,y)} a(x,y) f(y) dy.$$
Then for $\la >\delta^{-3}$
$$\|\cJ_\la\|_{L^2\to L^2}\lc \delta^{-1/2} \lambda^{-d/2}$$
where the implicit constants depend on $C_1$ in \eqref{uniflower}
and of a finite number of the constants in \eqref{unifupper}
($|\alpha|\le 10 d$ suffices).
\end{lemmasub}

\begin{proof}
We let $\delta_1= M^{-1}\delta$ where $M$ is very  large
in comparison to the constants in the assumptions (but independent of $\delta$ and $\la$).
By a partition of unity we may assume that the symbol $a$ is supported in the smaller cube
$Q\equiv Q_{\delta_1}(P^o)$.

We make affine  changes of variables  in $x$ and $y$ separately
which do not affect the assumptions, so that we may assume that
$P^o=O$, $\Psi_{x'y'}(O)=I_{d-1}$ (the $(d-1)\times(d-1)$ identity
matrix), and also $\Psi_{x'y_d}(O)=0$, $\Psi_{y'x_d}(O)=0$.

Then $\det (\Psi_{xy})= \Psi_{x_dy_d}+ O(\delta_1)$ in $Q$ and
thus $|\Psi_{x_d y_d}|\approx|\det(\Psi_{xy})|\approx \delta$. We
shall use orthogonality arguments based on the following
inequalities, valid for $(x,y)\in Q$ and $(x,z)$ in $Q$:
\begin{equation}\label{Psixprimediff}
|\Psi_{x'}(x,y)-\Psi_{x'}(x,z)|\geq |y'-z'| \quad\text{ if }
|y'-z'|\geq C_0\delta_1|y_d-z_d|,
\end{equation}
for a large constant $C_0$
and
\begin{equation}\label{Psixddiff}
|\Psi_{x_d}(x,y)-\Psi_{x_d}(x,z)|\geq C^{-1}\delta|y_d-z_d| \quad\text{ if }
|y'-z'|\leq
c_0\delta|y_d-z_d|;
\end{equation}
for a small constant $c_0$ but  $\delta_1$ is  so small that
$c_0\delta\gg C_0\delta_1$. Similar bounds hold for the phase
$\Psi^*(x,y):=\Psi(y,x)$. Inequality \eqref{Psixprimediff} follows
by a straightforward expansion about the origin, and it is crucial
that we use $\Psi_{x'y_d}(O)=0$. For \eqref{Psixddiff} we use of
course the lower bound on $\Psi_{x_dy_d}$.

We now decompose the amplitude into functions supported on rectangles
$R_m\times R_n$ (with $(m,n)\in \bbZ^d\times \bbZ^d$)
where both $R_m$ and $R_n$ have dimensions about
$\lambda^{-1/2}\times \dots\times \la^{-1/2}\times \la^{-1/2}\delta^{-1/2}$.
Let $\chi\in C^\infty_0(\bbR)$ so that $\chi $ is supported in $(-5/4,5/4)$ and $\sum_{j\in\bbZ} \chi(s-j)=1$ for all $s\in \bbR$.
Define for $(m,n)\in \bbZ^d\times \bbZ^d$
\begin{multline*}a_{m,n}(x,y)= \\a(x,y)
\chi(\la^{1/2}\delta^{1/2}x_d-m_d)
\chi(\la^{1/2}\delta^{1/2}y_d-n_d)
\prod_{i=1}^{d-1}
\chi(\la^{1/2}x_i-m_i))
\prod_{i=1}^{d-1}\chi(\la^{1/2}y_i-n_i)
\end{multline*}
and let $T_{mn}$ be defined as $\cJ_\la$ but with $a$ replaced by $a_{mn}$. Then
$\cJ_\la=\sum_{m,n} T_{mn}$.
We observe by simply using the support properties of the symbol and Schur's lemma that
\begin{equation}\label{Tind}
\|T_{mn}\|_{2-2}\lc \lambda^{-d/2}\delta^{-1/2};
\end{equation}
moreover by disjointness of symbols
\begin{equation}\begin{aligned}
\label{trivorth}
&T_{pq}^* T_{mn}=0 \quad&\text{ if $|p-m|\geq 4$},
\\
&T_{pq} T_{mn}^*=0\quad&\text{  if $|q-n|\geq 4$}.
\end{aligned}\end{equation}
In order to use the Cotlar-Stein orthogonality lemma it suffices to show
 that
\begin{align} \label{TstarT}
\|T_{pq}^* T_{mn}\|_{2-2}&\lc
\lambda^{-d}\delta^{-1} |n-n'|^{-N},
\text{ for $|m-p|\leq C$ and $|n-q|\geq C$,}
\\
 \label{TTstar}
\|T_{pq} T_{mn}^*\|_{2-2}&\lc
\lambda^{-d}\delta^{-1} |m-p|^{-N},
\text{ for $|n-q|\leq C$ and $|m-p|\geq C$.}
\end{align}
Let $H_{mnpq}(y,z)$ be the Schwartz kernel of $T_{pq}^*
T_{mn}$. By integration by parts we obtain the pointwise  bounds
$$
|H_{mnpq}(y,z)|\lc \lambda^{-d/2}\delta^{-1/2}\Big(
\frac{
\tfrac{|y_d-z_d|}{|y'-z'|}+ \sqrt{\la\delta}}
{\la|y'-z'|}\Big)^N\quad
 \text{ if }
|y'-z'|\geq
C_0\delta_1|y_d-z_d|,
$$
and
$$
|H_{mnpq}(y,z)|
\lc \lambda^{-d/2}\delta^{-1/2}
\Big(
\frac{\delta^{-2} +\sqrt{\la\delta}}{\la|y_d-z_d|}\Big)^N\quad
\text{ if }
|y'-z'|\leq
c_0\delta|y_d-z_d|.
$$
and since $C_0\delta_1\ll c_o\delta$ all relevant situations are
covered. In the first case we have $|y'-z'|\approx \la^{-1/2}
|m'-p'|$ and $|m'-p'|\gc |m_d-p_d|$, and in the second case we
have $|y_d-z_d|\approx \lambda^{-1/2}\delta^{-1/2} |m_d-p_d|$ and
$|m'-p'|\lc |m_d-p_d|$. By taking the support properties in $y$
and $z$ into account we can use Schur's Lemma to see that
\begin{align*}
&\sup_y\int|H_{mnpq}(y,z)| dz
+
\sup_z\int|H_{mnpq}(y,z)| dy
\\
&\lc \lambda^{-d}\delta^{-1} \begin{cases}
\Big(\frac{\delta^{-1}+\lambda^{1/2}\delta^{1/2}}
{\lambda^{1/2}|m'-p'|}\Big)^{N} \text{ if $|m'-p'|\geq c
|m_d-p_d|$},
\\
\Big(
\frac{\delta^{-3/2}+\lambda^{1/2} \delta}
{\lambda^{1/2}|m_d-p_d|}
\Big)^N \text{ if $|m'-p'|\leq C
|m_d-p_d|$.}
\end{cases}
\end{align*}
Our restriction  $\delta\ge \la^{-1/3}$ implies the desired bound
\eqref{TstarT}
for the operator norm of
$T_{pq}^* T_{mn}$. The operators
$T_{pq} T_{mn}^*$ are handled analogously.
\end{proof}

We now gather some facts that
are useful for $L^4$ estimates related to the Carleson-Sj\"olin theorem.
Define
\begin{align}
U_1\equiv U_1(s,\ts,t,\tt)&= s-t+\ts-\tt
\\
U_2\equiv U_2(s,\ts,t,\tt)&= (s-t)^2+(\ts-\tt)^2-2(t-\tt)(\ts-\tt)
\end{align}
and
\begin{equation}
\widetilde U_i(s,\ts,t,\tt)= U_i(\ts,s,\tt, t);
\end{equation}
moreover
\begin{align}
V_i(s,\ts,t,\tt)&=U_i(t,\tt, s,\ts)
\\
\widetilde V_i(s,\ts,t,\tt)&=\widetilde U_i(t,\tt, s,\ts)
\end{align}
(observe that  $U_1=\widetilde U_1=-V_1=-\widetilde V_1$).
The following calculus lemma is directly  taken from p. 63 in \cite{MuS}:

\begin{lemmasub}\label{calclemma}
 Let $A=(A_1,A_2)$ be an $\bbR^2$-valued function
 of class $C^{4}$,
defined on an interval.
Suppose that $M, M'>0$ and that
$ 2^{-M-1}\leq s-\ts\leq 2^{-M+1}$,
$ 2^{-M'-1}\leq t-\tt\leq 2^{-M'+1}$.
Let
\begin{equation}B(s,\ts,t,\tt)
=A(s)+A(\ts)-A(t)-A(\tt).\end{equation}
Then
(i)
\begin{equation}
\label{upBbound}
|B(s,\ts, t,\tt)|\leq C \min\{|U_1|+|U_2|, |V_1|+|V_2|\}
\end{equation}

(ii)
If also
\begin{equation}\label{a1a2}
|A_1'(s)A_2''(s)-A_2'(s) A_1''(s)|\geq c_1 \end{equation}
then
there is a uniform lower bound
\begin{equation}
\label{lowBbounddiffm}
|B(s,\ts, t,\tt)|\geq c \max \{|U_1|+|U_2|, |V_1|+|V_2|\} \quad\text { if  $|M-M'|>10$.}
\end{equation}

(iii) There are constants $c>0$,  $C_1>1$ so that if $M=M'$ then  the estimate
\begin{equation}\label{lowBboundsamem}
|B(s,\ts, t,\tt)|\geq c \max \{|U_1|+|U_2|, |V_1|+|V_2|\}
\end{equation}
holds in each of the following cases:
\begin{align}
&|s-t|+|\ts-\tt|\leq C_1^{-1}2^{-M}
\label{casea}
\\
\label{caseb}
\text{ or }
&|s-t|\geq C_1 2^{-M}
\\
\label{casec}
\text{ or }&|\ts-\tt|\geq C_1 2^{-M}.
\end{align}

(iv)
There is a constant $C_2>1$ so that the following holds.
Suppose  that either $M<M'-20$ or $M=M'$ and
$|s-t|+|\ts-\tt|\geq C_2 2^{-M} $. Suppose that
in addition $|U_1(s,\ts,t,\tt)|
\leq 2^{-M-10}$. Then
\begin{equation}
\label{U2MllMprime}
 |U_2(s,\ts,t,\tt)|\geq \frac 12 (\ts-\tt)^2\geq 2^{-2M-20}.
\end{equation}

(v) Suppose  $M=M'$ and let $\delta\leq 2^{-M-4}$. Suppose that
$|U_1(s,\ts,t,\tt)|\leq \delta/4$ and suppose that
$|s-t|+|\ts-\tt|\geq \delta$.
Then
$|s-t|\approx|\ts-\tt|$ and
\begin{equation}\label{U2MequalMprime}
 |U_2(s,\ts,t,\tt)|\geq 2^{-M-1}| \ts-\tt| \geq c 2^{-M}\delta
\end{equation}
\end{lemmasub}

\subsection{Proof of the orthogonality lemmata}
For the proofs of
Lemmata \ref{firstortho},
\ref{1.5ortho},
and
\ref{secondortho},
we shall need to analyze the expression
\begin{equation}
\label{pairingone}
\inn
{\fS^{m\mu b\tb}_{a\nu}
F^{m\mu b\tb}_{a\nu}}
{\fS^{m'\mu' b'\tb'}_{a'\nu'}
F^{m'\mu' b'\tb'}_{a'\nu'}}
\end{equation}
for the three cases $|m-m'|\geq 20$, then $m=m'$ and $|\mu-\mu'|\gg C$ and
finally $m=m'$, $\mu=\mu'$ and $|b-\tb|+|b'-\tb'|\gg C2^{-m-l/2}$.
We shall apply the lower bounds  of Lemma \ref{calclemma} (ii)
(with $M=l/2-m$, $M'=l/2-m'$)
 to the functions
$$y_1\mapsto A(y_1):=\phi_x(x, y_1, g(x,y_1))$$
and the upper bounds of  Lemma \ref{calclemma}
 to higher $x$ derivatives of $\phi$, evaluated at $y=(y_1,g(x,y_1))$.
The crucial  Carleson-Sj\"olin type condition \eqref{a1a2}
holds, as by a straightforward calculation
using \eqref{phixdy} and \eqref{mixedydterm}
$$
A_1'(y_1)A_2''(y_1)-A_2'(y_1) A_1''(y_1)=
\phi_{x_1y_1}\phi_{x_2y_1y_1} \Big|_{(x, y_1, g(x,y_1))} +o(1)
$$
where the $o(1)$ terms vanish at $O$ and the main terms are bounded below by
the  curvature
 condition \eqref{curvcond} (in the reduced form \eqref{curvatPreduced}).

Now we use the notation $s=y_1$, $t=z_1$, $\ts=\ty_1$, $\tt=\tz_1$, and
$U_1\equiv U_1(y_1,\ty_1,z_1,\tz_1)$ etc. Then
we have
\begin{multline*}
|\partial_x^\alpha\phi(x, y_1, g(x,y_1))+
\partial_x^\alpha\phi(x, \ty_1, g(x,\ty_1))
-\partial_x^\alpha \phi(x, z_1, g(x,z_1))-
\partial_x^\alpha\phi(x, \tz_1, g(x,\tz_1))|
\\
\leq C_{\alpha}\min
\{|U_1|+|U_2|, |V_1|+|V_2|\}
\end{multline*}
and
 in the cases
(i) $m'<m-10$  and (ii) $m=m'$ and one of
 \eqref{casea},  \eqref{caseb}, \eqref{casec}
we also get the lower bounds
\begin{multline*}
|\phi_x(x, y_1, g(x,y_1))+\phi_x(x, \ty_1, g(x,\ty_1))
-\phi_x(x, z_1, g(x,z_1))-\phi_x(x, \tz_1, g(x,\tz_1))|
\\
\geq
c\max
\{|U_1|+|U_2|, |V_1|+|V_2|\}
\end{multline*}
In the four term expressions that occur in the phases when writing
out \eqref{pairingone}  we have to replace $g(x,y_1)$ with $y_2$
etc. and then take into account that  $(x,y)$ belongs to $\supp
\zeta_l$; this introduces error terms of size $O(2^{-l})$.

Assuming that all points
$(x,y)$, $(x,z)$, $(x,\ty)$, $(x,\tz)$ belong to the support of
$\zeta_l$ then we obtain
\begin{multline}\label{4termupper}
|\partial_x^\alpha\phi(x, y)+
\partial_x^\alpha\phi(x, \ty)
-\partial_x^\alpha \phi(x, z)-
\partial_x^\alpha\phi(x, \tz)|
\leq C_{\alpha}\big(\min
\{|U_1|+|U_2|, |V_1|+|V_2|\}+  2^{-l}\big);
\end{multline}
moreover in the cases described above we also get the lower bound
\begin{multline}\label{4termlower}
|\phi_x(x, y)+\phi_x(x, \ty)
-\phi_x(x, z)-\phi_x(x, \tz)|
\geq
c\max
\{|U_1|+|U_2|, |V_1|+|V_2|\}
- C 2^{-l}.
\end{multline}
In order to further bound below the right hand side of
\eqref{4termlower} we shall use the statements in part (iv) and
(v) of Lemma \ref{calclemma}. It will turn out that in all the
described cases  $|U_1|+|U_2| \gg 2^{-l}$ so that the error terms
in \eqref{4termupper} and \eqref{4termlower} will not affect the
integrations by parts. This is an important point of the proof,
and many of our decompositions have been made with this goal in
mind.

Finally, before we discuss the proofs of the lemmata
we note that in all cases we may  assume
that $F^{m\mu b\tb}_{a\nu}$ is supported on a set of measure $2^{-2l-3  m}$;
simply  replace  $F^{m\mu b\tb}_{a\nu}$ with
$\widetilde\Omega^{m\mu b\tb}_{a\nu} F^{m\mu b\tb}_{a\nu}$
where $\widetilde \Omega^{m\mu b\tb}_{a\nu}
\Omega^{m\mu b\tb}_{a\nu}=
\Omega^{m\mu b\tb}_{a\nu}$ and
 $\widetilde \Omega^{m\mu b\tb}_{a\nu} $ has support properties similar to
 $\Omega^{m\mu b\tb}_{a\nu}$. Thus
\begin{equation}\label{CaSchw}
\|F^{m\mu b\tb}_{a\nu}\|_1
\lc 2^{-3l/2-m} \|F^{m\mu b\tb}_{a\nu}\|_2.
\end{equation}

\begin{proof}[Proof of Lemma \ref{firstortho}]

We square  the right hand side of \eqref{m-orthogonality} and see
 that we need to analyze
\eqref{pairingone} with
$|m-m'|\geq 20$. By symmetry we may assume that $m'<m-20$ (\textit{i.e. }
$2^{-l/2+m'}\ll 2^{-l/2+m}$).


We also apply part (iv) of Lemma \ref{calclemma}
which tells us that
in the present situation  $|U_1|+|U_2|\geq c 2^{-2M}\equiv c 2^{-l+2m}$.
We integrate by parts and observe that if derivatives hit the symbols involved we get a factor of $2^l$ with each derivative.

The size of the support of $u^{m\mu}_{a\nu}$ is $O(2^{-m-3l/2})$.
Consequently, after integrating by parts $2N$ times,   we obtain the bound
\begin{align*}
\big| &\inn
{\fS^{m\mu b\tb}_{a\nu}
F^{m\mu b\tb}_{a\nu}}
{\fS^{m'\mu' b'\tb'}_{a'\nu'}
F^{m'\mu' b'\tb'}_{a'\nu'}}
\big|
\\ \quad &\leq C_N (\la 2^{-2l+m})^{-2N} \meas\big( \supp u^{m\mu}_{a\nu}
\cap u^{m'\mu'}_{a'\nu'}\big) \big\|F^{m\mu b\tb}_{a\nu}\big\|_1
\big\|F^{m'\mu' b'\tb'}_{a'\nu'}\big\|_1
\\ \quad&\leq C_N' (\la 2^{-2l+m})^{-2N} 2^{-m-3l/2}
\big\|F^{m\mu b\tb}_{a\nu}\big\|_1
\big\|F^{m'\mu' b'\tb'}_{a'\nu'}\big\|_1.
\end{align*}
Now by the $T^*T$ argument using also the  Cauchy-Schwarz
inequality (for the terms with $|m-m'|\leq 20$) the expression on
the left hand side of  \eqref{m-orthogonality} is dominated by
$I+\sqrt{II}$ where $I$ is the first term on the right
 hand side of  \eqref{m-orthogonality} and
\begin{equation}\label{IIdef}
II= \sum_{0<m'<m-20<l/2}\sum_{\mu\in 2^{-l/2+m}\bbZ \atop {\mu'\in
2^{-l/2+m'}\bbZ} } \sum_{(b,\tb)\in \cP_\mu^m \atop (b',\tb')\in
\cP_{\mu'}^{m'} } \sum_{a,\nu, a',\nu'} \big|\inn {\fS^{m\mu
b\tb}_{a\nu} F^{m\mu b\tb}_{a\nu}} {\fS^{m'\mu' b'\tb'}_{a'\nu'}
F^{m'\mu' b'\tb'}_{a'\nu'}} \big|.
\end{equation}
We also  observe that for each fixed $m,m',\mu,\mu'$
the sums in  $(a,\nu)$ and $(a',\nu')$
are taken over index sets
 of cardinality $O(2^{3l/2+m})$ and $O(2^{3l/2+m'})$,
respectively. Moreover
for each fixed $m,\mu$ the sums in $(b,\tb)$ are over a set of
cardinality $2^{4m}$, and for each fixed $m',\mu'$ the sums in
$(b',\tb')$ are over a set of cardinality $2^{4m'}$. Finally for
each fixed $m$ the sums in $\mu$ and $\mu'$ are over sets of
cardinalities $O(2^{l/2 -m})$ and $O(2^{l/2 -m'})$, respectively.
Taking these restrictions into account
we continue with  straightforward estimation using just the
Cauchy-Schwarz inequality in the various parameters which gives an additional factor of
$2^{3l/2+m/2+m'/2} 2^{2m+2m'}2^{l/2-m/2-m'/2}$.
We  thus bound $|II|$ by
\begin{align*}
& C_N  \sum_{m,m'\atop{0<m'<m-20}}(\la 2^{-2l+m})^{-2N} 2^{-m-3l/2}
\sum_{\mu\in 2^{-l/2+m}\bbZ
\atop {\mu\in 2^{-l/2+m'}\bbZ} }
\sum_{(b,\tb)\in \cP_\mu^m
\atop (b',\tb')\in \cP_{\mu'}^{m'} }
\sum_{(a,\nu)
\atop{(a',\nu')}}
\big\|F^{m\mu b\tb}_{a\nu}\big\|_1
\big\|F^{m'\mu' b'\tb'}_{a'\nu'}\big\|_1
\\&
\leq C_N'  2^{l/2}(\la 2^{-2l})^{-2N}\sum_{0<m<l/2} 2^{(3-2N)m}
\sum_{\mu\in 2^{-l/2+m}\bbZ}
\sum_{(b,\tb)\in \cP_\mu^m}\sum_{(a,\nu)}
\big\|F^{m\mu b\tb}_{a\nu}\big\|_1^2.
\end{align*}
The assertion \eqref{m-orthogonality}
follows  if we  choose $N$ large in the previous estimate and apply
 \eqref{CaSchw}.
 \end{proof}

\begin{proof}[Proof of  Lemma \ref{1.5ortho}]
Now $m$ is fixed and we need to bound \eqref{pairingone} for $m=m'$ and
$|\mu-\mu'|\geq C 2^{m-l/2}$
for some large but absolute constant $C$.
We argue as in the proof of Lemma \ref{firstortho}, but now use
Lemma \ref{calclemma}, part (iii), \eqref{caseb} or \eqref{casec},
 with $M=l/2-m$. Thus the lower bound in
\eqref{4termlower} holds and also the upper bound in
\eqref{4termupper}. For the lower bounds we have $|U_1|+|U_2|\geq
c 2^{-2M}\approx 2^{-l+2m}$. Thus we get for $|\mu-\mu'| \geq C
2^{m-l/2}$,
\begin{align*}
\big| &\inn
{\fS^{m\mu b\tb}_{a\nu}
F^{m\mu b\tb}_{a\nu}}
{\fS^{m\mu' b'\tb'}_{a'\nu'}
F^{m\mu' b'\tb'}_{a'\nu'}}
\big|
\leq C_N (\la 2^{-2l+m})^{-2N} 2^{-m-3l/2}
\big\|F^{m\mu b\tb}_{a\nu}\big\|_1
\big\|F^{m\mu' b'\tb'}_{a'\nu'}\big\|_1.
\end{align*}
From here we proceed  as
 in the proof of Lemma \ref{firstortho}; we use the
Cauchy-Schwarz inequality in the parameters $\mu$, $\mu'$, $(b,\tb)$,
$(b',\tb')$,  and $\nu$,  and then  \eqref{CaSchw}.
\end{proof}

{\it  Remark.} One could also use
 Fourier transform arguments (with respect to $x$) as in
 the proof of Proposition  \ref{orthom=0} below.

\medskip

\begin{proof}[Proof of  Lemma \ref{secondortho}]
We have now $m$, $\mu$, $a$ and $\nu$ fixed, and we are required
to estimate $\|\sum_{b,\tb\in \cP_\mu^m} \fS^{m\mu
b\tb}_{a\nu}F^{m\mu b\tb}_{a\nu}\|_2$. The relevant $(b,\tb)$ is
such that $|b-\mu|\leq C2^{-l/2+m}$ and $|\tb-\mu|\leq
C2^{-l/2+m}$.

We split the family of pairs $\cP^m_\mu$ into a bounded set of
subfamilies
$\cP^m_{\mu,i}$
 with the property that
for any two pairs $(b,\tb)$, $(b',\tb')$ in one such
$\cP^m_{\mu,i}$ we have both $|b-b'|\leq c2^{-l/2+m}$ and
$|\tb-\tb'|\leq c2^{-l/2+m}$ for a small  constant $c$.

This time we need to  analyze \eqref{pairingone} with $m=m'$,
$\mu=\mu'$ and
$(b,\tb)\in \cP^m_{\mu,i}$,  $(b',\tb')\in \cP^m_{\mu,i}$.
We may use integration by parts
since
by the definition of $\cP^m_{\mu,i}$
we are in the situation of
part (iii), \eqref{casea} of Lemma \ref{calclemma},
with
$M=l/2-m$. The lower bound $|U_1|+|U_2|\geq 2^{-M}\delta$ in
\eqref{U2MequalMprime} applies with
\begin{equation}\label{bpairsdist}
\delta\approx |b-b'|+|\tb-\tb'|\ge C_4 2^{-m-l/2},\end{equation}
for some large $C_4$.
Thus in this case
\begin{align*}
\big|&\inn
{\fS^{m\mu b\tb}_{a\nu}
F^{m\mu b\tb}_{a\nu}}
{\fS^{m\mu b'\tb'}_{a'\nu'}
F^{m\mu b'\tb'}_{a'\nu'}}\big|
\\
&\quad\lc 2^{-m-3l/2}(\lambda 2^{-m-3l/2}
(|b-b'|+|\tb-\tb'|))^{-2N}
\big\|F^{m\mu b\tb}_{a\nu}\big\|_1
\big\|F^{m\mu b'\tb'}_{a'\nu'}\big\|_1
\\
\end{align*}
for $(b,\tb)\in \cP^m_{\mu,i}$,  $(b',\tb')\in \cP^m_{\mu,i}$
satisfying \eqref{bpairsdist}.
By a straightforward convolution inequality
\begin{multline*}
\Big\| \sum_{(b,\tb)\in \cP_\mu^m}
\fS^{m\mu b\tb}_{a\nu}F^{m\mu b\tb}_{a\nu}\Big\|_2 \lc\\
\Big(\sum_{(b,\tb)\in \cP_\mu^m} \big\|
\fS^{m\mu b\tb}_{a\nu}F^{m\mu b\tb}_{a\nu}\big\|_2^2\Big)^{1/2}
+
2^{-m/2-3l/4} (\la2^{-m-3l/2})^{-N}
\Big(\sum_{b,\tb}\|F^{m\mu b\tb}_{a\nu}\|_1^2\Big)^{1/2},
\end{multline*}
and \eqref{CaSchw} is used to  obtain the
 desired conclusion.
\end{proof}

\subsection{Proof of Proposition \ref{fixedfSlemma}}

This is to be deduced from Lemma \ref{basicestimate}.
%
%
We change variables in the integral defining
$\fS^{m\mu b\tb}_{a\nu}$ to
$$y_2=g(\tau_\nu^\mu a,y_1)+\sigma
, \quad \ty_2=g(\tau_\nu^\mu a,y_1)+\tsigma,
$$
where then integrations over $\sigma$, $\tsigma$ are extended over intervals of length $O(2^{-l})$.

We then have
\begin{equation}\label{reproffS}
\fS^{m\mu b\tb}_{a\nu}F(x)=
\big(\eta(2^{l}(\inn{\tfrac{\nabla_xg(a,\mu)}{|\nabla_xg(a,\mu)|}}{x-\tau^\mu_\nu
a}))\big)^2 \iint_{|\si|,|\tsi|\lc 2^{-l}} \cT^{m\mu b\tb}_{a\nu,
\sigma\tsi} [H_{\si\tsi}F]d\sigma d\tsigma
\end{equation}
where
\begin{equation}\label{Hsitsi}
H_{\sigma\tsigma} F(y_1,\ty_1)=
\beta_0(C^{-1}2^{l/2}\sigma)
\beta_0(C^{-1}2^{l/2}\tsigma)
F(y_1, g(\tau_\nu^\mu a,y_1)+\sigma, \ty_1,
g(\tau_\nu^\mu a,\ty_1)+\tsigma).
\end{equation}
The oscillatory integral operators
$\cT^{m\mu b\tb}_{a\nu, \sigma\tsi}$ in \eqref{reproffS}
act on functions $h$ of the variables
$(y_1,\ty_1)$ and are  defined by
\begin{equation}
\cT^{m\mu b\tb}_{a\nu,\sigma\tsi}h(x)=
\iint
\cA(x, y_1,\ty_1;\si,\tsi)
e^{i\la \Psi(x,y_1,\ty_1;\sigma,\tsi)}
h(y_1,\ty_1)
dy_1d\ty_1\,
\end{equation}
with
\begin{multline}
\cA(x, y_1,\ty_1;\si,\tsi)= (u^{m}_{a}(x))^2
\chi_1(2^{l/2-m-5}(y_1-\ty_1))
\zeta_l(x,y_1, g(\tau_\nu^\mu a,y_1)+\sigma)
\\ \times   \zeta_l(x,\ty_1,
g(\tau_\nu^\mu a,y_1)+\tsigma)
\,\Omega^{m b\tb}_{a}(y_1,g(\tau_\nu^\mu a,y_1)+\sigma,\ty_1,
g(\tau_\nu^\mu a,y_1)+\tsigma)
\label{amplitude}
\end{multline}
and
\begin{equation*}
\begin{gathered}
\Psi(x,y_1,\ty_1;\sigma,\tsi)=
\phi(x,y_1,g(\tau_\nu^\mu a,y_1)+\sigma)+
\phi(x,\ty_1,g(\tau_\nu^\mu a,\ty_1)+\tsi).
\end{gathered}
\end{equation*}

By the Cauchy-Schwarz inequality
\begin{equation}\label{afterC-S}
\|\fS^{m\mu b\tb}_{a\nu} F\|_2\lc 2^{-l}
\Big(\iint_{|\sigma|, |\tsi|\lc 2^{-l}}\big\|
\cT^{m\mu b\tb}_{a\nu, \sigma\tsi}
[H_{\si\tsi}F]\big\|_2^2d\sigma d\tsigma\Big)^{1/2}
\end{equation}
One now verifies that Lemma  \ref{basicestimate} with
$\delta\approx 2^{m-l/2}$ can be applied
to the operators
$\cT^{m\mu b\tb}_{a\nu, \sigma\tsi}$
so that the right hand side of
\eqref{afterC-S} is estimated by a constant times
\begin{equation*}
2^{-m/2-3l/4}
\Big(\iint_{|\sigma|, |\tsi|\lc 2^{-l}}\big\|
H_{\si\tsi}F\big\|_2^2d\sigma d\tsigma\Big)^{1/2}
\lc
2^{-m/2-3l/4} \|F\|_2.
\end{equation*}
\qed

\section{Estimation of  $\fB^0(f,f)$}\label{B0bound}
This case is handled rather analogously to the case $m>0$,
except instead of using Proposition \ref{fixedfSlemma} we reduce directly to
the Carleson-Sj\"olin theorem.

We shall set $\fS^{\mu b\tb}_{a,\nu}:=\fS^{0\mu b\tb}_{a,\nu}$,
$\cP_\mu:=\cP^0_\mu$, moreover $f^{\mu b}_{a\nu}(y):=\omega^{0\mu
b}_{a\nu}(y)f(y)$ ({\it cf.} \eqref{tensorF}), and define the
expressions $\cL f^{\mu b}_{a\nu}(y)$ and $\cE f^{\mu
b}_{a\nu}(y)$ by setting $m=0$ in  \eqref{Lfrequencyloc} and
\eqref{Efrequencyloc}. (Note that $\cP_\mu$ contains boundedly
many elements for each $\mu$.)

Note that Lemma \ref{yerrorterms} remains valid for $m=0$
so that $\fS^{0\mu b\tb}_{a,\nu}$ essentially acts on
 $\cL f^{\mu b}_{a\nu}\otimes
\cL f^{\mu \tb}_{a\nu}$.
Various orthogonality arguments will be  used for the proof of
\begin{proposition}\label{orthom=0} For $r\ge 2$ and
$N\in\mathbb{N}$,
\begin{multline} \label{ortho0eq}
\Big\|\sum_{\mu\in 2^{-l/2}\bbZ} \sum_{(b,\tb)\in \cP_\mu} \sum_{a,\nu}
\fS^{0\mu b\tb}_{a,\nu}
(f^{\mu b}_{a\nu}\otimes
f^{\mu\tb}_{a\nu})
\Big\|_r \lc\\
\Big( \sum_{a,\nu}
\sum_{\mu, b} \big\|
\fS^{0\mu b\tb}_{a,\nu}
(\cL f^{\mu b}_{a\nu}\otimes
\cL f^{\mu\tb}_{a\nu})
\big\|_r^{r'}\Big)^{1/r'}
+ C_N
\la^3(\la 2^{-2l})^{-N/r}
\sup_{a,\nu,\mu,b}
 \|f^{\mu b}_{a\nu}\|_{1}^{2}.
\end{multline}
\end{proposition}

We combine this with an application of the Carleson-Sj\"olin theorem
which will give
\begin{proposition}\label{CSappl} For $1\le p<4$, $r= 3p'/2$,
\begin{equation}
\big\|\fS^{0\mu b\tb}_{a,\nu}
(f\otimes g)\big\|_r \lc 2^{-2l/p'} \la^{-4/(3p')}
\|f\|_{p} \|g\|_{p};
\end{equation}
moreover
\begin{equation}
\big\|\fS^{0\mu b\tb}_{a,\nu}
(f\otimes g)\big\|_2 \lc 2^{-3l/2}
\la^{-1}  (\log\la)^{1/2}
\|f\|_{4} \|g\|_{4}.
\end{equation}
\end{proposition}

The error term in  \eqref{ortho0eq} is easily bounded by the
right hand side of \eqref{Bpqestimate} or \eqref{B0estimate} given that
$2^{l}\lc \la^{1/3}$.
For the main term in  \eqref{ortho0eq}
we  apply
Proposition \ref{CSappl} with $f\otimes g =
(\cL f^{\mu b}_{a\nu}\otimes
\cL f^{\mu\tb}_{a\nu})$ and put the result into
\eqref{ortho0eq}; this yields
\begin{align}
\Big( \sum_{a,\nu}
\sum_{\mu, b} \big\|
\fS^{0\mu b\tb}_{a,\nu}
(\cL f^{\mu b}_{a\nu}\otimes
\cL f^{\mu\tb}_{a\nu})
\big\|_r^{r'}\Big)^{1/r'}
\lc &A_p(\la, l) \Big(\sum_{a,\nu}\sum_{\mu,b}
\big\|\cL f^{\mu b}_{a\nu}\big\|_p^{2r'}\Big)^{1/r'}
\notag
\\
\label{lpLpsq}
\lc &A_p(\la, l) \Big(\sum_{a,\nu}\sum_{\mu,b}\big\|\cL f^{\mu b}_{a\nu}\big\|_p^{p}\Big)^{2/p},
\end{align}
where $A_p(\la, l)=
2^{-2l/p'} \la^{-4/(3p')}$ if $p>4$ and
$A_4(\la, l)= 2^{-3l/2} \la^{-1}  (\log\la)^{1/2} $.
In the last displayed inequality we have used that if
$r=3p'/2$ then  $2r'\ge p$ holds for $p\le 4$.
 The desired estimate for $\fB_0(f,f)$ then follows from an application of Lemma \ref{lpLplemma} to \eqref{lpLpsq}.

\begin{proof}[Proof of Proposition \ref{orthom=0}]

The $u^0_a$ are supported on cubes $Q_a$
with diameter $\approx 2^{-l/2}$, centered at $a$,
which are essentially disjoint
(so that $\sum_a\chi_{Q_a}(x)\le C$).
Thus
\begin{multline*}
\Big\|
\sum_{\mu\in 2^{-l/2}\bbZ} \sum_{(b,\tb)\in \cP_\mu} \sum_{a,\nu}
\fS^{0\mu b\tb}_{a,\nu}
(f^{\mu b}_{a\nu}\otimes
f^{\mu\tb}_{a\nu})
\Big\|_r    \\
\lc
\Big(\sum_a \Big\|\chi_{Q_a}\sum_{\mu\in 2^{-l/2}\bbZ} \sum_{(b,\tb)\in \cP_\mu}
 \sum_{\nu}
\fS^{0\mu b\tb}_{a,\nu}
(f^{\mu b}_{a\nu}\otimes
f^{\mu\tb}_{a\nu})
\Big\|_r^r\Big)^{1/r}.
\end{multline*}
Let $\eta_0\in C_{0}^{\infty}(\mathbb{R}^2)$ be such that
$\eta_0(s)=1$ if $|s|\le 1$. Let $W_{a\mu} $ be the convolution
operator on functions in $\bbR^2$ which has Fourier multiplier
$$w_{a\mu}(\xi)=\eta_0(C^{-1} \la^{-1}2^{l/2}(\xi-2\lambda \nabla_x\phi(a,\mu, g(a,\mu))));$$
here $C$ is chosen so large that
$|\lambda^{-1}\xi-\nabla_x\phi(x,y)-\nabla_x\phi(x,\ty)|\ge
2^{-l/2}$ whenever $w_{a\mu}(\xi)=0$ and $(x,y,\ty)$ is in the
convex hull of the support of $u^m_a(x) \om^{mb}_a(y_1)
\om^{m\tb}_a(\ty_1)$ for all $(b,\tb)\in \cP_\mu$.

In view of this property we obtain by the inversion formula for
the Fourier transform and a straightforward integration by parts
argument that
\begin{eqnarray}
\begin{aligned}
|(I-W_{a\mu})\fS^{0\mu b\tb}_{a,\nu} (f^{\mu b}_{a\nu}\otimes
f^{\mu\tb}_{a\nu})(x)|& \leq \|(1-w_{a\mu})(\fS^{0\mu
b\tb}_{a,\nu} (f^{\mu b}_{a\nu}\otimes
f^{\mu\tb}_{a\nu}))\:\widehat{\;}\:\|_1
\\&\le C_N\lambda^2 2^{-2l} (\lambda
2^{-3l/2})^{-N} \|f^{\mu b}_{a\nu}\|_1 \|f^{\mu\tb}_{a\nu}\|_1
\end{aligned}
\end{eqnarray}
for all $N\in\mathbb{N}$, uniformly in $x\in\mathbb{R}^2$. From
this the contribution
$$
\Big(\sum_a \Big\|\chi_{Q_a}
\sum_{\mu\in 2^{-l/2}\bbZ} \sum_{(b,\tb)\in \cP_\mu}
 \sum_{\nu}
(I-W_{a\mu})
\fS^{0\mu b\tb}_{a,\nu}
(f^{\mu b}_{a\nu}\otimes
f^{\mu\tb}_{a\nu})
\Big\|_r^r\Big)^{1/r}
$$
can be estimated by the error term in \eqref{ortho0eq} in a
straightforward way (we use that there are $O(2^{l})$ relevant
$a$'s, $O(2^{l/2})$ relevant $\mu$'s, and for fixed $a,\mu$ there
are $O(2^{l/2})$ relevant $\mu$'s and $O(1)$ relevant  $b$'s).

For the main term we use the orthogonality properties of the operators $W_{a\mu}$
(with respect to $\mu$ when  $a$ is
fixed) and then the essential disjoint support of the functions $u^{0\mu}_{a,\nu}$ (when $a,\mu$ are fixed). We obtain
\begin{multline}\Big(\sum_a \Big\|\chi_{Q_a}
\sum_{\mu\in 2^{-l/2}\bbZ} \sum_{(b,\tb)\in \cP_\mu}
 \sum_{\nu}
W_{a\mu} \fS^{0\mu b\tb}_{a,\nu} (f^{\mu b}_{a\nu}\otimes
f^{\mu\tb}_{a\nu}) \Big\|_r^r\Big)^{1/r}\\ \lc \Big(\sum_a\Big(
\sum_{\mu\in 2^{-l/2}\bbZ} \sum_{(b,\tb)\in \cP_\mu}
 \sum_{\nu}
\Big\|
\fS^{0\mu b\tb}_{a,\nu}
(f^{\mu b}_{a\nu}\otimes
f^{\mu\tb}_{a\nu})
\Big\|_r^{r'}\Big)^{r/r'}\Big)^{1/r}.
\end{multline}

Finally using Lemma \ref{yerrorterms} we can replace
$(f^{\mu b}_{a\nu}\otimes
f^{\mu\tb}_{a\nu})$ by
$(\cL f^{\mu b}_{a\nu}\otimes
\cL f^{\mu\tb}_{a\nu})$ as the other terms just contribute to the error term in \eqref{ortho0eq}.
\end{proof}

\begin{proof}[Proof of Proposition \ref{CSappl}]
Define
\begin{equation}\label{Tdef}
T^{\mu b}_{a\nu} f(x) = u^{0\mu}_{b\nu}(x)\int \zeta_l(x,y)
\om^{0\mu b}_{a\nu}(y) e^{i\la \phi(x,y) } f(y) dy
\end{equation}
We  dispose of
 the diagonal  cutoff function $\chi_0(2^{l/2 -5}(y_1-\ty_1))$ in the definition of $\fS^{0\mu b\tb}_{a\nu}$ by expanding $\chi_0$
in a Fourier series and obtain
\begin{equation*}
\fS^{0\mu b\tb}_{a\nu}(f\otimes g) =\sum_{k\in \bbZ} c_k
T^{\mu b}_{a\nu} f_k
(x)  T^{\mu b}_{a\nu} g_k(x)
\end{equation*}
where $f_k(y)= f(y)e^{ik 2^{l/2-5}y_1}$, $g_k(y)= g(y)e^{-ik
2^{l/2-5}y_1}$ and $|c_k|\le C_N 2^{-N|k|}$ for all
$N\in\mathbb{N}$. Then
$$
\big\|\fS^{0\mu b\tb}_{a,\nu}
(f\otimes g)\big\|_r \lc \sup_k \big\|T^{\mu b}_{a\nu} f_k\big\|_{2r}
\big\|T^{\mu \tb}_{a\nu} g_k\|_{2r}.
$$
We now change variables $y_2=g(a,b)+\sigma$ in \eqref{Tdef} and in view of the support assumption the $\sigma$ integration is extended over an interval of length $\le C2^{-l}$.
The phase $\Psi^\sigma(x,y_1)=\Phi(x,y_1, g(a,b)+\sigma)$
is a phase satisfying the assumptions of the Carleson-Sj\"olin theorem
with bounds uniform in the parameters.
Thus if we set
\begin{equation*}\label{Tsigdef}
T^{\mu b,\sigma}_{a\nu} h(x) = u^{0\mu}_{b\nu}(x)\int
\zeta_l(x,y_1,g(a,b)+\sigma) \om^{0\mu b}_{a\nu}(y_1,
g(a,b)+\sigma) e^{i\la \phi(x,y_1,g(a,b)+\sigma) } h(y_1) dy_1
\end{equation*}
we obtain
with  $2r=3p'$, $p<4$
\begin{align*}
\big\|T^{\mu b}_{a\nu} f_k\big\|_{2r}
&\le \int_{|\sigma|\le C2^{-l}}
\big\|T^{\mu b,\sigma}_{a\nu} [f_k(\cdot, g(a,b)+\sigma)]\big\|_{2r}
 d\sigma
\\&\lc \la^{-2/2r}
\int_{|\sigma|\le C2^{-l}}
\big\|f_k(\cdot, g(a,b)+\sigma)\big\|_{p} d\sigma
\\&
\lc 2^{-l/p'} \la^{-3/p'}\|f_k\|_p,
\end{align*}
and of course $\|f_k\|_p=\|f\|_p$. We argue similarly in the case $p=2r=4$
but then one gets an additional factor $(\log \la)^{1/4}$ in the bound.
\end{proof}

\section{Appendix I: \\ A sharpening of an $L^p$ improving
inequality  for averages on curves}\label{appI}

Consider the translation invariant averaging operator
\begin{equation}
\cA f(x)=\int \chi(s) f(x-\gamma(s)) ds
\end{equation}
for the curve
\begin{equation}
\gamma(s)=\big(s, \frac{s^2}2, \frac{s^3}6\big),
\end{equation}
where $\chi$ denotes a cutoff function to a neighborhood of $0$.
The sharp $L^p\to L^q$ estimates are known and due to Oberlin
\cite{Ob}, in fact $\cA$ maps $L^p\to L^q$ if and only if
$(1/p,1/q)$ belongs to the trapezoid with corners $(0,0)$,
$(1,1)$, $(1/2, 1/3)$ and $(2/3,1/2)$. However the
 critical  $L^2\to L^{3}$ and $L^{3/2}\to L^2$ estimates
can be improved if one uses Lorentz spaces;  this improvement
does to the best of our knowledge not follow from the $T^*T$ method used in \cite{Ob}.

\begin{theorem}
$\cA$ maps $L^2(\bbR^3)$  to $L^{3,2}(\bbR^3)$ and
$L^{3/2,2}(\bbR^3)$ to $L^2(\bbR^3)$.
\end{theorem}


\begin{proof}
By duality it suffices to prove the $L^2\to L^{3,2}$ inequality.
By a standard reduction using Littlewood-Paley  theory it suffices
to prove that the operators $\cA_k$ defined by
$$
\cA_k f(x)=\int f(y) \iint
e^{i\sigma(x_2-y_2-\tfrac{(x_1-y_1)^2}{2})+
i\tau(x_3-y_3-\tfrac{(x_1-y_1)^3}{6})}
\chi_1(\tfrac{\sigma^2+\tau^2}{2^{2k}} )\chi(x-y) d\sigma d\tau
\,dy
$$
map  $L^2$ to $L^{3,2}$ boundedly (with norms uniformly in
$k\gg0$). Here $\chi_1\in C_0^{\infty}(\mathbb{R})$ is an
appropriate cutoff function supported away from $0$.
The reason for
the validity of this reduction is that
 $\cA_k = L_k\cA_k L_k+E_k$
where $L_k$ are Littlewood-Paley operators localizing frequencies to annuli of width $C2^{k}$ and the errors $E_k$ satisfy
$\|E_k\|_{L^p\to L^q}=O(2^{-k})$.
Then assuming that
\begin{equation}\label{Akpiece}
\sup_k \|\cA_k\|_{L^2\to L^{3/2}}\leq A,
\end{equation}
we obtain
\begin{align*}
&\Big\|\sum_{k>0} L_k \cA_k L_k f\Big\|_{L^{3,2}}
\lc
\Big\|\Big(\sum_{k>0}|\cA_k L_k f|^2\Big)^{1/2}\Big\|_{L^{3,2}}
\\
&\lc\Big(\sum_{k>0}\|\cA_k L_k f\|_{L^{3,2}}^2\Big)^{1/2}
\lc A
\Big(\sum_{k>0}\| L_k f\|_{L^{2}}^2\Big)^{1/2} \lc A\|f\|_2
\end{align*}

For the first inequality we used Littlewood-Paley theory, and for the
second one we used a Minkowski-type inequality which amounts to the imbedding $\ell^2(L^{3/2}) \subset L^{3,2} (\ell^2)$ which can be seen using
the equivalence
$\|u^2\|_{L^{p/2,q/2}}\approx \|u\|_{L^{p,q}}^2$ and
the triangle inequality in the Lorentz-space $L^{3/2,1}$.

We now turn our attention to the operators $\cA_k$ and the proof of
\eqref{Akpiece}.
 We are fortunate as our inequality involves the space $L^2$ at
least on the function side and
one can reduce matters to the estimation of
 an oscillatory integral operator
\begin{equation} \label{Tla}
T_\la f(x) = \int e^{i\la \Phi(x,y)}\chi(x,y) f(y) dy
\end{equation}
mapping $L^2$ to $L^{3,2}$ with norm $O(\lambda^{-d/3})$. Here
$\la\approx 2^k$ and the phase is given by
\begin{equation}\label{phase}
\Phi(x,y) = y_2(x_2+\tfrac{(x_1-y_1)^2}{2})
+y_3(x_3+\tfrac{(x_1-y_1)^3}{6}),
\end{equation}
where  $|y_3| \ge c>0$ in the support of $\chi$.
The reduction to the oscillatory integral operator involves
Plancherel's theorem (with respect to the $(y_2,y_3)$ variables),
a rescaling by $2^k$, and renaming $(\sigma,\tau)$ to $(y_2,y_3)$.

Now define
$$T_{\la,l}  f(x) = \int e^{i\la \Phi(x,y)}\chi(x,y)
\chi(2^l \det \Phi_{xy}) f(y) dy
$$
where the cutoff function localizes to the set where
$|\det \Phi_{xy}|\approx 2^{-l}$. Observe that
$$
-\det \Phi_{xy}=  y_2+y_3(x_1-y_1).
$$
Define $\widetilde T_\la$ by a similar cutoff to the region where
$|\det \Phi_{xy}|\leq \la^{-1/3}$ and choose the cutoff functions so that
$T_\la=\sum_{2^l<\la^{1/3}}T_{\la, l} +\widetilde T_\la$.

We have the usual $L^2$ bounds
$\|T_{\la,l}\|_{L^2\to L^2} \lc 2^{l/2} \la^{-3/2}$ for $2^l<\la^{1/3}$
and $\|\widetilde T_\la\|_{L^2\to L^2} \lc\la^{-4/3}$.
Now it remains to show that for $r<4$, $s=3r'$
\begin{align}\label{rsineq1}
\|T_{\la,l}\|_{L^r\to L^s} &\lc 2^{-l/r'} \la^{-3/r}, \quad 2^l<\la^{1/3}
\\
\label{rsineq2}
\|\widetilde T_\la\|_{L^2\to L^2} &\lc \la^{-3/r-1/3r'}
\end{align}

This implies that $T_\la$ maps $L^{p,1}$ to $L^{q,\infty}$ for
$1\leq p<5/2$, $q=3p'/2$. By another real interpolation we deduce
that $T_\la$ maps in fact $L^{p,a}$ to $L^{q,a}$ for any $a>0$,
and choosing $a=p$ and $p=2$ yields the $L^2\to L^{3,2}$
inequality for $T_\la$.

The proofs of  \eqref{rsineq1},  \eqref{rsineq2}  follow the ideas
in Theorem \ref{main2d}, however as in \S4 we may directly reduce
matters to B. Barcel\'o's  restriction theorem for cones
(\cite{B}). In fact let  for $w\in \Bbb R^3$, $v\in \Bbb R^2$,
$v_2\approx 1$
$$\Psi(w,v)= v_2\big(-w_1v_1^2+w_2v_1+w_3),$$
and define $S_\la f(w)=\int e^{i\la \Psi(w,v)} \chi(w,v)f(v) dv$.
Then
\begin{equation} \label{conerestr}
\|S_\la f\|_s\lc \la^{-3/s} \|f\|_r, \quad s=3r', \, r<4.
\end{equation}
Indeed from the restriction theorem we get \eqref{conerestr} for
cutoff's of the product form $\chi(w,v)=a(w) b(v)$, and by using a
Fourier series expansion of $\chi$ we can reduce to this case.

In order to prove \eqref{rsineq1} we split $f=\sum_{\ka\in \bbZ}
f_\ka$ where the function $f_\ka$ is supported in $\{y: |-y_1+y_2/y_3-\ka
2^{-l}|\le 2^{-l}\}$,  and then we observe that $T_{\la, l} f_\kappa$
is supported where $|x_1-\ka 2^{-l}|\le C 2^{-l}$. It thus suffices to
show \eqref{rsineq1} for $f$ replaced with $f_\kappa$. Fix
$\kappa$ and let $\alpha_\kappa=\kappa 2^{-l}$. We change
variables in $y$ by setting $y_2=-y_3(\alpha_\ka-y_1+\sigma)$
(where  $|\sigma|\lc 2^{-l}$) and then set $v_1=y_1$, $v_2=y_3$.
After a short computation we obtain  $$ \Phi(x;y_1,
-y_3(\alpha_\ka-y_1+\sigma), y_3)-y_1^3y_3/3=
\Psi(h^\sigma(x),y_1,y_3)\equiv \Psi(h^\sigma(x),v)
$$
with
$$
h^\sigma(x)=(x_1+\sigma-\alpha_\ka,x_2-(\sigma-\alpha_\ka)x_1-x_1^2/2,
x_3+(\sigma-\alpha_\ka)(x_2+x_1^2/2)+x_1^3/6);
$$ clearly as a nonlinear shear $h^\sigma$ defines a global diffeomorphism.

Now
$$
T_{\la,l} f_\ka(x)= \int_{|\sigma|\lc 2^{-l}}
\int e^{i\la \Psi(h^\sigma(x),v)}
\chi^\sigma(x,v) f_{\ka,\sigma}(v) dv
$$
with
$\chi^\sigma(x,v) f_{\ka,\sigma}(v) =
\chi(x,v_1,-v_2(\alpha_\ka-v_1+\sigma),v_2)
f_\ka(v_1,-v_2(\alpha_\ka-v_1+\sigma),v_2) e^{i\la v_1^3v_2/3}$.
Thus from
\eqref {conerestr} we get
\begin{align*}
\|T_{\la,l}f_\ka\|_{L^s(\Bbb R^3)} \lc\int_{|\sigma|\lc 2^{-l}}
\|f_{\ka,\sigma}\|_{L^r(\Bbb R^2)} d\sigma
\lc 2^{-l/r'}
\Big(\int
\|f_{\ka,\sigma}\|_{L^r(\Bbb R^2)}^r d\sigma\Big)^{1/r} \lc
 2^{-l/r'} \|f_\ka\|_{L^r(\Bbb R^3)}
\end{align*}
which implies \eqref{rsineq1}. Inequality \eqref{rsineq2} is
proved in a similar way. \end{proof}


\section{Appendix II:
\\
 On bilinear versions of an adjoint restriction theorem for the circle}
\label{appII}

Let $d\sigma$ denote arclength measure on the circle and for $f\in
L^2(\bbS^1)$ consider the family of (restricted) extension
operators given by $\cE_\la f(\theta) = \widehat{fd\sigma}(\la
\theta)$, where $\theta\in\mathbb{S}^1$. In \cite{BBC1} Barcel\'o,
Carbery and one of the current authors proved the sharp bilinear
inequality
\begin{equation}\label{BBCineq}
\int  |\cE_\la f(\theta)  \cE_\la g(\theta) | d\sigma(\theta)
\leq C \la^{-5/6}\|f\|_{L^2(\bbS^1)} \|g\|_{L^2(\bbS^1)}
\end{equation}
valid under the separation conditions

\begin{equation}\label{sep}
\supp (f)\cap\supp(g)=\emptyset \;\text{ and } \;\supp (\widetilde
f)\cap\supp(g)=\emptyset;
\end{equation}
here $\widetilde f(\theta):=f(-\theta)$. The initially complicated
proof of  this inequality in \cite{BBC1} has since been simplified
in \cite{BBC2}. Here we generalize and simplify further this
result by interpreting the separation condition \eqref{BBCineq} as
a condition on the associated canonical relation, and deduce the
bilinear estimates directly from known linear estimates of the
type \eqref{L2Tlestimate}, \eqref{L2Tlestimatelim} (proved already in
\cite{PS} for $d=1$). Several multilinear extensions of this
argument are possible (in the spirit of \cite{BBC2})
 but we shall not pursue this here.

Consider the oscillatory integral operator defined on
functions in $L^2(\bbR)$ by
$$T_\la f(x)=\int e^{i\la\phi(x,y)} \chi(x,y) f(y) dy$$
where $\chi$ is smooth and compactly supported. We assume that the
canonical relation $\{(x,\phi_x, y, -\phi_y)\}$ is a folding
canonical relation, \textit{i.e.} \eqref{piLfold} and
\eqref{piRfold} hold. We may clearly assume that $\phi_{xy}$   is
small  on the support of the amplitude.
%

\begin{proposition}\label{bilosc}
Suppose that $\chi$ is supported on a set of diameter at most
$\delta_0$  and suppose that
 $\phi_{xyy}\neq 0$ and $\phi_{xxy}\neq 0$ on the support of $\chi$.

Let $0<\delta<\delta_0$ and suppose that
\begin{equation}\label{separation}\dist (\supp f, \supp g) \geq \delta>0.
\end{equation}
Then, if $\delta_0$ is  sufficiently small then, for  large $\la>0$,

(i)
\begin{equation}\label{claim1}
\big \| T_\la f \, T_\la g\big\|_1 \leq C_\delta \la^{-5/6} \|f\|_2\|g\|_2.
\end{equation}

(ii) Moreover

\begin{equation}\label{claim2}
\big \| T_\la f \, T_\la g\big\|_{6/5} \leq C_\delta \la^{-5/6}
\big[ \|f\|_2\|g\|_3+ \|f\|_3 \|g\|_2\big] .
\end{equation}
\end{proposition}

{\it Remark:} The inequality
 \eqref{BBCineq} (valid assuming \eqref{sep})
can be deduced by applying
the proposition with the phase $\phi(x,y)=\cos (x-y)$.

\begin{proof} [Proof of Proposition \ref{bilosc}]
We may assume that $\la\gg \delta^{-1}$.
A better inequality follows immediately by the standard  $L^2$ estimates
if we assume that $\phi_{xy}\neq 0$; thus  we shall assume that $\phi_{xy} $ vanishes somewhere and after a straightforward reduction
(using
suitable localizations)  we may assume that on the support of the
relevant cutoff function
$\phi_{xy}=0 \iff  x=h(y)$
where $h$ is invertible and $|h'|$ is bounded above and below
(actually, as in \cite{PS} one can reduce to $h(y)=y$, by a change of
variable in
$y$).

Let $\chi_0\in C^\infty_0(\bbR)$ be supported in $(-1,1)$ and equal to one in $(-1/2,1/2)$, and let
 $\chi_1(t)=\chi(2t)-\chi(t)$, $\chi_l(t)=\chi_1(2^l t)$.
Consider the operators given by
$T_{\la, l}$ and $\widetilde T_\la$ defined in and following
 \eqref{Tlaldef}.
Notice that the kernel $K_{\la, l}(x,y)$ is supported where
$ |x-h(y)|\approx  2^{-l}$, and
and similarly the kernel for $S_\la$ is supported where $|x-h(y)|\lc \la^{-1/3}$.

Now write
\begin{equation*}
T_\la f(x) T_\la g(x)=
\Big( \widetilde T_{\la} f(x)+\sum_{2^l\leq \la^{1/3}} T_{\la, l} f(x)\Big)
\Big( \widetilde T_{\la} g(x)+\sum_{2^m\leq \la^{1/3}} T_{\la, m} g(x)\Big).
\end{equation*}
We  use the  separation assumption on the supports of $f$ and $g$.
If $y\in \supp f$ and $z\in \supp g$ then the conditions
$ |x-h(y)|\approx  2^{-l}$,   and
$ |x-h(z)|\approx  2^{-m}$,  can hold  {\it simultaneously} only if
either $l\leq \ell_0$, or $m\leq \ell_0$,
for some fixed $\ell_0=\ell_0(\delta)$.
Thus
$$T_{\la, l} f(x)  T_{\la, m} g(x)=0 \quad \text{ if } l\geq \ell_0, \text{ and } m\geq \ell_0.$$
Similarly
$\widetilde T_\la f(x) \widetilde T_{\la} g(x)=0$, and
$\widetilde T_\la f(x) T_{\la, l} g(x)=0 $, $l\leq \ell_0$.
Therefore,
$$
T_\la f T_\la g= I_\la(f,g)+ II_\la (f,g) +III_\la (f,g)
$$
where
\begin{align*}
I_\la(f,g) &=\Big( \widetilde T_{\la} f+\sum_{2^{\ell_0}<2^l\leq \la^{1/3}}
T_{\la, l} f\Big)
\Big( \sum_{m\leq \ell_0(\delta)} T_{\la, m} g\Big),
\\
II_\la(f,g) &=\Big( \sum_{l\leq \ell_0(\delta)} T_{\la, l} f\Big)
\Big( \widetilde T_{\la} g+\sum_{2^{\ell_0}<2^m\leq \la^{1/3}}
T_{\la, l} g\Big),
\\
III_\la(f,g) &=\Big( \sum_{l\leq \ell_0(\delta)} T_{\la, l} f\Big)
\Big( \sum_{m\leq \ell_0(\delta)} T_{\la, l} g\Big).
\end{align*}
We now use the  Cauchy-Schwarz inequality and apply
the standard $L^2$  bounds
\eqref{L2Tlestimate}, \eqref{L2Tlestimatelim} (with $d=1$). We obtain
\begin{align*}
\|I_\la(f,g)\|_1&\leq
\Big(
\| \widetilde T_{\la} f\|_2 +
\sum_{2^{\ell_0}<2^l\leq \la^{1/3}}
\|T_{\la, l} f\|_2\Big) \sum_{m\leq \ell_0(\delta)} \|T_{\la, m} g\|_2
\\
&\leq C(\delta)  \la ^{-1/2}  (\la^{-1/3}+
\sum_{2^{\ell_0}<2^l\leq \la^{1/3}} 2^{l/2} \la^{-1/2})
\|f\|_2\|g\|_2
\leq C_1(\delta) \la^{-5/6} \|f\|_2\|g\|_2 .
\end{align*}

Similarly
one proves
the inequality
$\|II_\la(f,g)\|_1\leq C_2(\delta) \la^{-5/6} \|f\|_2\|g\|_2$,
and also
the bound
  $\|III_\la(f,g)\|_1\leq C_3(\delta)\la^{-1} \|f\|_2\|g\|_2$.
This shows \eqref{claim1}.

For \eqref {claim2} we use  the endpoint
$L^3\to L^3$ inequality  in \cite{GS3} which says that
$$\Big\|\sum_{2^{\ell_0}<2^l\leq \la^{1/3}}
T_{\la, l} f\Big\|_3\leq C_\delta \la^{-1/3} \|f\|_3
$$
and also $\|S_\la\|_{L^3\to L^3}=O(\la^{-1/3})$.
Since $L^2\cdot L^3\subset L^{6/5}$
we apply H\"older's inequality to obtain
\begin{align*}
\|I_\la(f,g)\|_{6/5} &\leq C_\delta \la^{-5/6} \|f\|_3\|g\|_2,
\\
\|II_\la(f,g)\|_{6/5} &\leq C_\delta \la^{-5/6} \|f\|_2 \|g\|_3,
\\
\|III_\la(f,g)\|_{6/5} &\leq C_\delta \la^{-5/6} \|f\|_2\|g\|_2,
 \end{align*}
and thus \eqref{claim2}.
\end{proof}
\medskip

\end{document}